\documentclass{amsart}
\usepackage{amsfonts}

\usepackage{graphicx}
\usepackage{epstopdf}

\input tcilatex

\begin{document}
\title[Genealogy of branching processes]{On the genealogy and coalescence times of Bienaym\'{e}-Galton-Watson
branching processes}
\author{Nicolas Grosjean, Thierry Huillet$^{*}$}
\address{Laboratoire de Physique Th\'{e}orique et Mod\'{e}lisation \\
CNRS-UMR 8089 et Universit\'{e} de Cergy-Pontoise, 2 Avenue Adolphe Chauvin,
95302, Cergy-Pontoise, FRANCE, E-mail: Nicolas.Grosjean@u-cergy.fr,
Thierry.Huillet@u-cergy.fr}
\maketitle

\begin{abstract}
Coalescence processes have received a lot of attention in the context of
conditional branching processes with fixed population size and
non-overlapping generations. Here we focus on similar problems in the
context of the standard unconditional Bienaym\'{e}-Galton-Watson branching
processes, either (sub)-critical or supercritical. Using an analytical tool,
we derive the structure of some counting aspects of the ancestral genealogy
of such processes, including: the transition matrix of the ancestral count
process and an integral representation of various coalescence times
distributions, such as the time to most recent common ancestor of a random
sample of arbitrary size, including full size.

We illustrate our results on two important examples of branching mechanisms
displaying either finite or infinite reproduction mean, their main interest
being to offer a closed form expression for their probability generating
functions at all times. Large time behaviors are investigated.\newline



$^{*}$ corresponding author.
\end{abstract}

\section{\textbf{Introduction and outline of the results}}

Using a sampling without replacement formula, we develop a general analytic
approach allowing to understand some aspects of the ancestral-count of the
discrete-time-$t$ Bienaym\'{e}-Galton-Watson process with current population
size $N_{t}\left( n_{0}\right) $, started with $n_{0}$ founders. The states
of the ancestral process correspond to the numbers of common ancestors, when
moving backwards in time, of a given sub-sample of the current population $%
N_{t}\left( n_{0}\right) $. We consider the (sub-)critical and supercritical
cases. In particular, we derive:\newline

- the analytic expression of the backward-in-time block-counting transition
matrix of the ancestral process, evaluating the one-step backward
probability, on the event $\left\{ N_{t}\left( n_{0}\right) \geq i\right\} $%
, to move from state $i$ to state $j\leq i$, through ancestral merging. In
sharp contrast with similar concern for Fisher-Wright like constant
population size branching models, \cite{MS}, \cite{Hui2}, this
(lower-triangular) transition matrix is time-inhomogeneous and
sub-stochastic. In the Fisher-Wright setup, with a very rich development
starting from \cite{King}, the starting point model is a conditional
branching process, introduced in \cite{KMG1} and \cite{KMG2}, as a
population model with fixed population size, and non-overlapping generations.

- the probability distribution of $\tau _{i,j}^{\left( t\right) }\left(
n_{0}\right) $, the first time, starting from $i$ randomly-chosen
individuals at time $t$, that the block-counting ancestral process ever
enters state $j\leq i$, as measured from generation $0$. Our results
complement and generalize the ones obtained in \cite{Lam} and \cite{Le}
studying the coalescence times $\tau _{i,1}^{\left( t\right) }\left(
n_{0}\right) $ for a finite number $i$ of individuals sampled in the current
generation in the subcritical case.

- the joint probability of the event $\tau _{i,1}^{\left( t\right) }\left(
n_{0}\right) \geq k,N_{t}\left( n_{0}\right) =j$. It gives the distribution
of the time-to-most-recent-common ancestor (TMRCA) of $i$ randomly sampled
individuals on the event $N_{t}\left( n_{0}\right) =j\geq i$.

- the probability of the event $\tau _{i,1}^{\left( t\right) }\left(
n_{0}\right) =\infty $, either because the $i$\ randomly sampled individuals
do not belong to the offspring of a common founder among $n_{0}$\ or because
there are strictly less than $i$\ individuals alive at generation $t$.

- the conditional probability distribution of $\tau _{i,1}^{\left( t\right)
}\left( 1\right) $, given $n_{0}=1$ and $N_{t}\left( 1\right) =i$, namely $%
\mathbf{P}\left( \tau _{i,1}^{\left( t\right) }\left( 1\right) \geq k\mid
N_{t}\left( 1\right) =i\right) $. It is the TMRCA given the whole
single-founder population alive at $t$ with size $i$\ is being sampled.

- the distribution of the coalescence time for the whole population of the $%
t $-th generation conditioned on the event that it is not extinct, which,
specifically, is: $\mathbf{P}\left( \tau _{N_{t}\left( 1\right) ,1}^{\left(
t\right) }\left( 1\right) \geq k\mid N_{t}\left( 1\right) >0\right) $.

- the ratio of triple versus binary one-step merging probabilities.\newline

In principle, our results would lead to general statements on the
asymptotics of the TMRCA under various conditionings and for various general
branching mechanisms. However, it turns out to be computationally involved
with such a degree of generality. We rather illustrate the results on some
generic explicit examples covering a wide range of situations, without
aiming at exhaustivity.

Our approach giving an integral representation of the probabilities of
interest is indeed particularly well-suited when the
probability-generating-function of the population size $N_{t}\left(
n_{0}\right) $ is available in closed-form for all times: after warming up
with $b$-ary deterministic trees, we proceed by illustrating our results on
the Bienaym\'{e}-Galton-Watson branching process with geometric reproduction
law (with finite mean number of offspring per capita) and the
Bienaym\'{e}-Galton-Watson branching process with `Sibuya' reproduction law
(with infinite mean number of offsprings). Both reproduction laws share an
invariance under iteration property. They are particular incarnations of a
family of generalized linear-fractional models introduced in \cite{Sag} and
further studied in \cite{GH}.

The geometric Bienaym\'{e}-Galton-Watson branching process is of particular
interest because its contour process is known to be the standard (fair or
unfair) Harris random walk, \cite{Har}, \cite{Ben}. In particular, following
our alternative algebraic path, we recover various large-$t$-limiting
results (both qualitatively and quantitatively) on the distribution of $\tau
_{2,1}^{\left( t\right) }\left( 1\right) $ conditioned on $N_{t}\left(
1\right) \geq 2$ and of $\tau _{N_{t}\left( 1\right) ,1}^{\left( t\right)
}\left( 1\right) $ conditioned on $N_{t}\left( 1\right) >0$. These results
are mirrored by the literature in the (sub)-critical, \cite{Lam}, \cite{A}, 
\cite{Le}, and supercritical cases, \cite{A3}. Specifically, the coalescence
time, both for pairs of tips and for the whole population, occur in the
recent past (in the subcritical regime), in the distant past (in the
supercritical regime) or in-between (in the critical regime).

The Bienaym\'{e}-Galton-Watson branching process with heavy-tailed Sibuya
branching mechanism is a prototype of an extreme branching process studied
in \cite{A2} from the point of view of coalescence. We recover some results
derived therein, to which we add that, given $N_{t}\left( 1\right) >0$, $%
\tau _{N_{t}\left( 1\right) ,1}^{\left( t\right) }\left( 1\right) $ has a
limiting geometric distribution.\newline

Let us briefly discuss the position of this work with respect to former
works on similar topics. In \cite{Zub}, while assigning to each particle a
``genealogy'', limiting distributions of the (rescaled or not) distance to
the closest common ancestor of any two particles are obtained, in the
supercritical, subcritical and $\left( \alpha ,L\right) -$critical cases,
both in discrete and continuous time. Extensions to multitype branching
processes are also supplied. In \cite{Lam}, similar results are obtained
both in discrete-time and continuous space/time settings, dealing
essentially with the subcritical case. Corollary $1$ in that paper was a
source of inspiration to develop the present analytic approach, based on
sampling without replacement from the current population. In \cite{A3},
using purely probabilistic tools, similar results on the probability
distribution of the time to most recent common ancestors of any two randomly
sampled individuals at generation $t$ and its behavior as $t\rightarrow
\infty $ under various general conditions on the branching mechanism, either
supercritical, critical or subcritical, are derived. In \cite{A}, in the
same spirit and setting, limiting distributions of the (rescaled or not)
coalescence time for the whole population conditioned on the event that it
is not extinct are also obtained. In \cite{A2}, it is shown that for general
rapidly growing populations (with infinite mean number of offspring per
capita), coalescence occurs in the recent past rather than remote past, the
latter case being rather typical of supercritical with finite mean branching
mechanisms. Finally, in \cite{Le}, limiting distributions of the distance to
the closest common ancestor of any two particles or more are obtained in the
context of a continuous-time Bienaym\'{e}-Galton-Watson branching processes.
The latter work seems unaware of the work \cite{Zub}, covering part of this
case.

As stated before, our work produces new general integral formulas of the
probabilities of the events of interest in this ancestral-count context, in
terms of the branching mechanism, its iterates and derivatives. As such, it
is particularly well-suited when the probability-generating-function of the
population size is available in closed-form for all times and our examples
offer a narrower scope than in the above-cited works. For instance, for
rapidly growing populations, our explicit Sibuya example does not cover the
case where this probability generating function has a slowly varying
function factor as in \cite{A2}. Note that we can deal with discrete-time
Bienaym\'{e}-Galton-Watson branching processes with any number of founders.

\section{\textbf{Count genealogies of Bienaym\'{e}-Galton-Watson branching
processes}}

In order to fix the notations, we first revisit well-known facts on
Bienaym\'{e}-Galton-Watson branching processes. Consider a discrete-time-$t$
branching process whose reproduction law has probability mass $\mathbf{P}%
\left( M=m\right) =\pi _{m},$ $m\geq 0$ for the number $M$ of offspring per
capita, see \cite{Harris} and \cite{AN}. We avoid the trivial case $\pi
_{1}=1$. We let $\phi \left( z\right) =\mathbf{E}\left( z^{M}\right) $ be
the probability generating function (pgf) of $M$ and we assume $\phi \left(
1\right) =1$ (no finite-time explosion). Let $N_{t}\left( n_{0}\right) $ be
the number (possibly $0$) of individuals alive at generation $t\geq 1$,
given $N_{0}=n_{0}\geq 1.$ We have 
\begin{equation*}
\phi _{t}\left( z\right) :=\mathbf{E}\left( z^{N_{t}\left( 1\right) }\right)
=\phi ^{\circ t}\left( z\right) ,
\end{equation*}
with $\phi ^{\circ t}\left( z\right) $ the $t-$th composition of $\phi
\left( z\right) $ with itself. The Bienaym\'{e}-Galton-Watson process $%
N_{t}\left( 1\right) $ is a time-homogeneous Markov chain with denumerable
state-space $\Bbb{N}_{0}:=\left\{ 0,1,...\right\} $. Furthermore, 
\begin{equation*}
\mathbf{E}\left( z^{N_{t}\left( n_{0}\right) }\right) =\phi
_{t}^{n_{0}}\left( z\right) ,
\end{equation*}
because $N_{t}\left( n_{0}\right) =\sum_{m=1}^{n_{0}}N_{t}^{\left( m\right)
}\left( 1\right) $, with $N_{t}^{\left( m\right) }\left( 1\right) $ the
descendance of the $m$-th founder at time $t$, all mutually independent. We
shall let $\mu =\mathbf{E}\left( M\right) $ (if this quantity is finite) so
that $\mathbf{E}\left( N_{t}\left( 1\right) \right) =\mu ^{t}$ and $\mathbf{E%
}\left( N_{t}\left( n_{0}\right) \right) =n_{0}\mu ^{t}$. The
Bienaym\'{e}-Galton-Watson process $N_{t}\left( n_{0}\right) $ has a
time-homogeneous stochastic transition matrix $P^{\left( n_{0}\right) }$,
with entries $P_{i,j}\left( n_{0}\right) =\left[ z^{j}\right] \phi \left(
z\right) ^{i}=\mathbf{P}\left( N_{1}\left( i\right) =j\right) $ (with $%
\left[ z^{j}\right] \phi \left( z\right) ^{i}$\ denoting the $z^{j}$%
-coefficient of the pgf $\phi \left( z\right) ^{i}$) and initial condition $%
\mathbf{P}\left( N_{0}\left( n_{0}\right) =j\right) =\delta _{n_{0},j}$.

Recalling that $N_{t}\left( 1\right) =\left( \sum_{l=1}^{N_{t-1}\left(
1\right) }M_{l}\right) \mathbf{1}_{N_{t-1}\left( 1\right) >0}$, $t\geq 1$,
is the number of individuals alive at generation $t$ given $n_{0}=1,$
therefore $\phi _{t}\left( z\right) $ obeys the recurrence \cite{AN}, 
\begin{equation*}
\phi _{t}\left( z\right) =\phi \left( \phi _{t-1}\left( z\right) \right) 
\text{, }\phi _{0}\left( z\right) =z,
\end{equation*}
with $\mathbf{P}\left( N_{t}\left( 1\right) =n\right) :=\left[ z^{n}\right]
\phi _{t}\left( z\right) $, the $z^{n}-$coefficient of the power series $%
\phi _{t}\left( z\right) $.

Depending on $\mu <1$ or $\mu \geq 1$, the process is called subcritical or
(super-)critical with almost sure finite-time extinction in the sub-critical
and critical cases only.

\subsection{Count genealogy of the Bienaym\'{e}-Galton-Watson process $%
\left\{ N_{t}\left( n_{0}\right) \right\} :$\ Results}

The count genealogy process is a time-inhomogeneous integral-valued Markov
process, say $\left\{ \widehat{N}_{t}\left( n_{0}\right) \right\} $, whose
lower-triangular transition matrix is given by its $\left( i,j\right) -$%
entries $\widehat{P}_{i,j}^{\left( t\right) }\left( n_{0}\right) $, giving
the probability to move from state $i\geq 1$ to state $j,$ $1\leq j\leq i$,
from generation $t$ to generation $t-1$ backwards in time. So time $t$ of $%
\left\{ \widehat{N}_{t}\left( n_{0}\right) \right\} $ varies backward,
starting from any generation and ending up at generation $0$. This
probability is obtained by taking an $i$-sample without replacement from $%
N_{t}\left( n_{0}\right) $, tracing back its ancestry (the number of its
ancestors one generation before) and evaluating the probability that this
number is $j$. If the branching process has $n_{0}$ founders, with $\left(
M\right) _{i}:=M\left( M-1\right) ...\left( M-i+1\right) $, the falling
factorial of $M$, on the set $\left\{ N_{t}\left( n_{0}\right) \geq
i\right\} $, we get 
\begin{equation}
\begin{array}{l}
\widehat{P}_{i,j}^{\left( t\right) }\left( n_{0}\right) :=\mathbf{E}\left( 
\binom{N_{t-1}\left( n_{0}\right) }{j}\sum_{i_{1}+...+i_{j}=i}^{*}\binom{i}{%
i_{1}...i_{j}}\prod_{l=1}^{j}\frac{\left( M_{l}\right) _{i_{l}}}{\left(
\sum_{l^{\prime }=1}^{N_{t-1}\left( n_{0}\right) }M_{l^{\prime }}\right)
_{i_{l}}}\right) \\ 
=\sum_{n\geq j}\mathbf{P}\left( N_{t-1}\left( n_{0}\right) =n\right) \binom{n%
}{j}\sum_{i_{1}+...+i_{j}=i}^{*}\binom{i}{i_{1}...i_{j}}\mathbf{E}\left(
\prod_{l=1}^{j}\frac{\left( M_{l}\right) _{i_{l}}}{\left( \sum_{l^{\prime
}=1}^{n}M_{l^{\prime }}\right) _{i_{l}}}\right) .
\end{array}
\label{1a}
\end{equation}
In Eq. (\ref{1a}), the star-sum is to indicate $i_{1},...,i_{j}\geq 1$.%
\newline

\emph{Remark 1:} Clearly, $\sum_{j=1}^{i}\widehat{P}_{i,j}^{\left( t\right)
}\left( n_{0}\right) =\mathbf{P}\left( N_{t}\left( n_{0}\right) \geq
i\right) $ because to take an $i$-sample without replacement from $%
N_{t}\left( n_{0}\right) $ requires $N_{t}\left( n_{0}\right) \geq i$. The
matrix $\widehat{P}_{i,j}^{\left( t\right) }\left( n_{0}\right) $, as
defined in (\ref{1a}), is sub-stochastic and does not as such define a
proper Markov process. To define one such process, one could complete the
state-space $\Bbb{N}:=\left\{ 1,2,...\right\} $ of $\left\{ \widehat{N}%
_{t}\left( n_{0}\right) \right\} $ by adding a coffin-state $\left\{
\partial \right\} $ in which one enters from $i$ with probability $\mathbf{P}%
\left( N_{t}\left( n_{0}\right) <i\right) $. The augmented transition matrix
is now stochastic with added state $\left\{ \partial \right\} $ absorbing
and state-space $\Bbb{N}\cup \left\{ \partial \right\} $, opening the way to
quasi-stationarity questions while conditioning on the event $N_{t}\left(
n_{0}\right) \geq i$. Note that starting from state $\left\{ 1\right\} $,
the only accessible states are $\left\{ 1,\partial \right\} $. $\Diamond $%
\newline

We give the following representation of $\widehat{P}_{i,j}^{\left( t\right)
}\left( n_{0}\right) $ in terms of the branching mechanism $\phi $, its
iterates and derivatives:

\begin{proposition}
\label{prop1} On the set $N_{t}\left( n_{0}\right) \geq i$, the
time-inhomogeneous transition matrix of the count genealogical process reads 
\begin{equation}
\begin{array}{l}
\widehat{P}_{i,j}^{\left( t\right) }\left( n_{0}\right) =\mathbf{P}\left( 
\widehat{N}_{t-1}\left( n_{0}\right) =j\mid \widehat{N}_{t}\left(
n_{0}\right) =i\right)  \\ 
=\frac{1}{j!\Gamma \left( i\right) }\sum_{i_{1}+...+i_{j}=i}^{*}\binom{i}{%
i_{1}...i_{j}}\int_{0}^{1}dz\left( 1-z\right) ^{i-1}\left[ \phi
_{t-1}^{n_{0}}\right] ^{\left( j\right) }\left( \phi \left( z\right) \right)
\prod_{l=1}^{j}\phi ^{\left( i_{l}\right) }\left( z\right) .
\end{array}
\label{1c}
\end{equation}
\newline
\end{proposition}

\begin{corollary}
\label{corol1} The following identity holds 
\begin{equation}
\begin{array}{l}
\mathbf{P}\left( N_{t}\left( n_{0}\right) \geq i\right)  \\ 
=\sum_{j=1}^{i}\frac{1}{j!\Gamma \left( i\right) }%
\sum_{i_{1}+...+i_{j}=i}^{*}\binom{i}{i_{1}...i_{j}}\int_{0}^{1}dz\left(
1-z\right) ^{i-1}\left[ \phi _{t-1}^{n_{0}}\right] ^{\left( j\right) }\left(
\phi \left( z\right) \right) \prod_{l=1}^{j}\phi ^{\left( i_{l}\right)
}\left( z\right) .
\end{array}
\label{1c2}
\end{equation}
\newline
\end{corollary}

\textbf{First hitting times of the backward count process.} With $j\leq i$,
let 
\begin{equation*}
\tau _{i,j}^{\left( t\right) }\left( n_{0}\right) =\sup \left\{ 0\leq s<t:%
\widehat{N}_{s}\left( n_{0}\right) =j\mid \widehat{N}_{t}\left( n_{0}\right)
=i\right\} ,
\end{equation*}
with the convention that $\tau _{i,j}^{\left( t\right) }\left( n_{0}\right)
=\infty $ if this set is empty. This set could be empty either as result of $%
N_{t}\left( n_{0}\right) <i$ (see Remark $1$) or because the $i$ sampled
individuals are the offspring of more than $j$ founders. $\tau
_{i,j}^{\left( t\right) }\left( n_{0}\right) $ is the last time $s$ that $%
\widehat{N}_{s}\left( n_{0}\right) $ enters state $j$ backward in time,
given $\widehat{N}_{t}\left( n_{0}\right) =i$, as measured from the
generation number $0$. In terms of the one-step transition probability of $%
\left\{ \widehat{N}_{s}\left( n_{0}\right) \right\} $, we have 
\begin{equation*}
\mathbf{P}\left( \tau _{i,j}^{\left( t\right) }\left( n_{0}\right)
=t-1\right) =\widehat{P}_{i,j}^{\left( t\right) }\left( n_{0}\right)
\end{equation*}
and more generally,\newline

\begin{theorem}
\label{theo1} For a given $n_{0}\geq 1$ and all $k\in \left\{
0,...,t-1\right\} $, the distribution of $\tau _{i,j}^{\left( t\right)
}\left( n_{0}\right) $\ is given by
\end{theorem}

\begin{equation}
\begin{array}{l}
\mathbf{P}\left( \infty >\tau _{i,j}^{\left( t\right) }\left( n_{0}\right)
\geq k\right) \\ 
=\frac{1}{j!\Gamma \left( i\right) }\sum_{i_{1}+...+i_{j}=i}^{*}\binom{i}{%
i_{1}...i_{j}}\int_{0}^{1}dz\left( 1-z\right) ^{i-1}\left[ \phi
_{k}^{n_{0}}\right] ^{\left( j\right) }\left( \phi _{t-k}\left( z\right)
\right) \prod_{l=1}^{j}\phi _{t-k}^{\left( i_{l}\right) }\left( z\right) .
\end{array}
\label{1d}
\end{equation}
\newline

\begin{corollary}
\label{corol2} For a given $n_{0}\geq 1$ and all $k\in \left\{
0,...,t-1\right\} $, the distribution of the TMRCA of $i$\ randomly sampled
individuals out of $N_{t}\left( n_{0}\right) $\ is given by 
\begin{equation}
\begin{array}{l}
\mathbf{P}\left( \infty >\tau _{i,1}^{\left( t\right) }\left( n_{0}\right)
\geq k\right)  \\ 
=\frac{1}{\Gamma \left( i\right) }\int_{0}^{1}dz\left( 1-z\right)
^{i-1}\left[ \phi _{t}^{n_{0}}\left( z\right) \right] ^{\left( 1\right) }%
\frac{\phi _{t-k}^{\left( i\right) }\left( z\right) }{\phi _{t-k}^{\left(
1\right) }\left( z\right) }.
\end{array}
\label{1de}
\end{equation}
While tracking jointly the population size $N_{t}\left( n_{0}\right) $\ at
time $t$, 
\begin{equation}
\begin{array}{l}
\mathbf{E}\left( \mathbf{1}_{\tau _{i,1}^{\left( t\right) }\left(
n_{0}\right) \geq k}u^{N_{t}\left( n_{0}\right) }\right)  \\ 
=\frac{u^{i}}{\Gamma \left( i\right) }\int_{0}^{1}dz\left( 1-z\right)
^{i-1}\left( \phi _{k}^{n_{0}}\right) ^{\left( 1\right) }\left( \phi
_{t-k}\left( uz\right) \right) \phi _{t-k}^{\left( i\right) }\left(
uz\right) ,
\end{array}
\label{1e}
\end{equation}
so that, with $j\geq i$, 
\begin{equation}
\begin{array}{l}
\mathbf{P}\left( \infty >\tau _{i,1}^{\left( t\right) }\left( n_{0}\right)
\geq k,N_{t}\left( n_{0}\right) =j\right)  \\ 
=\frac{1}{\Gamma \left( i\right) }\left[ u^{j-i}\right] \int_{0}^{1}dz\left(
1-z\right) ^{i-1}\left( \phi _{k}^{n_{0}}\right) ^{\left( 1\right) }\left(
\phi _{t-k}\left( uz\right) \right) \phi _{t-k}^{\left( i\right) }\left(
uz\right) .
\end{array}
\label{1f}
\end{equation}
\newline
\end{corollary}

\emph{Remark 2:}

Note that, putting $i=2$ and $k=t-1$ in the integral representation of $%
\mathbf{P}\left( \infty >\tau _{i,1}^{\left( t\right) }\left( n_{0}\right)
\geq k\right) $,

\begin{equation*}
\begin{array}{l}
\widehat{P}_{2,1}^{\left( t\right) }\left( n_{0}\right) =\mathbf{E}\left(
N_{t-1}\left( n_{0}\right) \frac{\left( M_{1}\right) _{2}}{\left(
\sum_{l^{\prime }=1}^{N_{t-1}\left( n_{0}\right) }M_{l^{\prime }}\right) _{2}%
}\right) =\int_{0}^{1}\left( 1-z\right) \left[ \phi _{t}\left( z\right)
^{n_{0}}\right] ^{\left( 1\right) }\frac{\phi ^{\left( 2\right) }\left(
z\right) }{\phi ^{\left( 1\right) }\left( z\right) }dz \\ 
=n_{0}\int_{0}^{1}\left( 1-z\right) \phi _{t}\left( z\right) ^{n_{0}-1}\phi
_{t}^{\left( 1\right) }\left( z\right) \frac{\phi ^{\left( 2\right) }\left(
z\right) }{\phi ^{\left( 1\right) }\left( z\right) }dz,
\end{array}
\end{equation*}
is the probability that two randomly chosen individuals at generation $t$
will merge immediately at generation $t-1$. The quantity $1-\widehat{P}%
_{2,1}^{\left( t\right) }\left( n_{0}\right) $ is the probability that they
will not merge. $\Diamond $\newline

Let $\tau _{2,1}^{\left( t\right) }\left( n_{0}\right) $ be the TMRCA of two
randomly chosen individuals at generation $t$ (the distance of their merging
time to the generation $0$), with the convention that $\tau _{2,1}^{\left(
t\right) }\left( n_{0}\right) =\infty $ if these two individuals either
belong to the descendance of two distinct founders among $n_{0}$ or if there
are less than two individuals alive at generation $t$. We have 
\begin{eqnarray*}
\tau _{2,1}^{\left( t\right) }\left( n_{0}\right) &=&t-1\text{ if the two
sampled particles merge in one step }t\rightarrow t-1 \\
&&\tau _{2,1}^{\left( t\right) }\left( n_{0}\right) \overset{d}{=}\tau
_{2,1}^{\left( t-1\right) }\left( n_{0}\right) \text{ if not.}
\end{eqnarray*}
Thus, with $\varphi _{t}\left( z\right) =\mathbf{E}\left( z^{\tau
_{2,1}^{\left( t\right) }\left( n_{0}\right) }\right) $, by first-step
analysis, 
\begin{equation*}
\varphi _{t}\left( z\right) =\widehat{P}_{2,1}^{\left( t\right) }\left(
n_{0}\right) z^{t-1}+\left( 1-\widehat{P}_{2,1}^{\left( t\right) }\left(
n_{0}\right) \right) \varphi _{t-1}\left( z\right) \text{, }\varphi
_{0}\left( z\right) :=0.
\end{equation*}
In the latter recurrence, $\varphi _{0}\left( z\right) $ has been set to $0$
because two randomly chosen individuals at generation $0$ have no common
ancestor, viz $\tau _{2,1}^{\left( 0\right) }\left( n_{0}\right) =\infty $.
This recurrence yields $\varphi _{t}\left( z\right) $ in principle. More
generally, let $\tau _{i,1}^{\left( t\right) }\left( n_{0}\right) $\ be the
TMRCA of $i$\ randomly sampled individuals out of $N_{t}\left( n_{0}\right) $%
\ (the distance of their merging time to the founders at generation $0$).
The latter recurrence can be generalized to $\varphi _{t,i}\left( z\right) :=%
\mathbf{E}\left( z^{\tau _{i,1}^{\left( t\right) }\left( n_{0}\right)
}\right) $ with $i>2$, giving the distribution of $\tau _{i,1}^{\left(
t\right) }\left( n_{0}\right) $ by first-step-analysis. Equation (\ref{1de})
is an alternative representation of this distribution. Furthermore, $\tau
_{i,1}^{\left( t\right) }\left( n_{0}\right) =\infty $\ if either these $i$\
individuals do not belong to the descendance of a common founder among $%
n_{0} $\ or if there are strictly less than $i$\ individuals alive at
generation $t $ and\newline

\begin{proposition}
\label{prop2}$\left( i\right) $ 
\begin{equation}
\mathbf{P}\left( \tau _{i,1}^{\left( t\right) }\left( n_{0}\right) <\infty
\right) =\frac{1}{\Gamma \left( i\right) }\int_{0}^{1}dz\left( 1-z\right)
^{i-1}\left[ \phi _{t}^{n_{0}}\left( z\right) \right] ^{\left( 1\right) }%
\frac{\phi _{t}^{\left( i\right) }\left( z\right) }{\phi _{t}^{\left(
1\right) }\left( z\right) }.  \label{finite}
\end{equation}
$\left( ii\right) $ The probability $\mathbf{P}\left( \tau _{i,1}^{\left(
t\right) }\left( n_{0}\right) <\infty \right) $ obeys the recurrence, 
\begin{equation}
\begin{array}{l}
\mathbf{P}\left( \tau _{i,1}^{\left( t\right) }\left( n_{0}\right) <\infty
\right) -\mathbf{P}\left( \tau _{i-1,1}^{\left( t\right) }\left(
n_{0}\right) <\infty \right)  \\ 
=-\mathbf{P}\left( N_{t}\left( n_{0}\right) =i-1\right) -n_{0}\left(
n_{0}-1\right) \int_{0}^{1}dz\left( 1-z\right) ^{i-1}\phi _{t}^{\left(
1\right) }\left( z\right) \phi _{t}^{\left( i-1\right) }\left( z\right) \phi
_{t}\left( z\right) ^{n_{0}-2}.
\end{array}
\label{ip}
\end{equation}
\newline
\end{proposition}

\emph{Remark 3:}

If $n_{0}=1$, the only possible reason why $\tau _{i,1}^{\left( t\right)
}\left( 1\right) =\infty $ is because there are strictly less than $i$
individuals alive at generation $t$. Therefore 
\begin{equation*}
\mathbf{P}\left( \tau _{i,1}^{\left( t\right) }\left( 1\right) <\infty
\right) =\mathbf{P}\left( N_{t}\left( 1\right) \geq i\right) .
\end{equation*}
Let us check this from (\ref{finite}) with $n_{0}=1:$ with $\phi
_{t}^{\left( i\right) }\left( z\right) =\sum_{n\geq i}\left( n\right)
_{i}z^{n-i}\mathbf{P}\left( N_{t}\left( 1\right) =n\right) $ and $%
\int_{0}^{1}dz\cdot z^{n-i}\left( 1-z\right) ^{i-1}=\Gamma \left( i\right)
/\left( n\right) _{i}$, 
\begin{eqnarray*}
\frac{1}{\Gamma \left( i\right) }\int_{0}^{1}dz\left( 1-z\right) ^{i-1}\phi
_{t}^{\left( i\right) }\left( z\right) &=&\frac{1}{\Gamma \left( i\right) }%
\sum_{n\geq i}\left( n\right) _{i}\mathbf{P}\left( N_{t}\left( 1\right)
=n\right) \int_{0}^{1}dz\cdot z^{n-i}\left( 1-z\right) ^{i-1} \\
&=&\mathbf{P}\left( N_{t}\left( 1\right) \geq i\right) .\text{ }\Diamond
\end{eqnarray*}
\newline

\begin{corollary}
\label{corol3} The probability that $\tau _{i,1}^{\left( t\right) }\left(
n_{0}\right) =\infty $\ as a result of these $i$\ individuals not belonging
to the descendance of a common founder among $n_{0}$\ only, is: 
\begin{equation}
\mathbf{P}\left( \tau _{i,1}^{\left( t\right) }\left( n_{0}\right) =\infty
,N_{t}\left( n_{0}\right) \geq i\right) =\mathbf{P}\left( N_{t}\left(
n_{0}\right) \geq i\right) -\mathbf{P}\left( \tau _{i,1}^{\left( t\right)
}\left( n_{0}\right) <\infty \right) .  \label{I1}
\end{equation}
In particular, if $i=2$, 
\begin{equation}
\begin{array}{l}
\mathbf{P}\left( \tau _{2,1}^{\left( t\right) }\left( n_{0}\right) =\infty
,N_{t}\left( n_{0}\right) \geq 2\right) =\mathbf{P}\left( N_{t}\left(
n_{0}\right) \geq 2\right) -\mathbf{P}\left( \tau _{2,1}^{\left( t\right)
}\left( n_{0}\right) <\infty \right)  \\ 
=n_{0}\left( n_{0}-1\right) \int_{0}^{1}dz\left( 1-z\right) \phi
_{t}^{\left( 1\right) }\left( z\right) ^{2}\phi _{t}\left( z\right)
^{n_{0}-2}.
\end{array}
\label{I2}
\end{equation}
\end{corollary}

Because what really makes sense is the probability distribution of $\tau
_{i,1}^{\left( t\right) }\left( n_{0}\right) $ given $\tau _{i,1}^{\left(
t\right) }\left( n_{0}\right) <\infty $, we can also state\newline

\begin{corollary}
\label{corol4} For each $n_{0}\geq 1$\ and $k\in \left\{ 0,...,t-1\right\} $%
\begin{equation}
\mathbf{P}\left( \tau _{i,1}^{\left( t\right) }\left( n_{0}\right) \geq
k\mid \tau _{i,1}^{\left( t\right) }\left( n_{0}\right) <\infty \right) =%
\frac{\int_{0}^{1}dz\left( 1-z\right) ^{i-1}\left[ \phi _{t}^{n_{0}}\left(
z\right) \right] ^{\left( 1\right) }\frac{\phi _{t-k}^{\left( i\right)
}\left( z\right) }{\phi _{t-k}^{\left( 1\right) }\left( z\right) }}{%
\int_{0}^{1}dz\left( 1-z\right) ^{i-1}\left[ \phi _{t}^{n_{0}}\left(
z\right) \right] ^{\left( 1\right) }\frac{\phi _{t}^{\left( i\right) }\left(
z\right) }{\phi _{t}^{\left( 1\right) }\left( z\right) }}.  \label{I3}
\end{equation}
When $n_{0}=1$ (a single founder), $\tau _{i,1}^{\left( t\right) }\left(
1\right) <\infty \Leftrightarrow N_{t}\left( 1\right) \geq i$. Furthermore,\ 
\begin{equation}
\mathbf{P}\left( \tau _{i,1}^{\left( t\right) }\left( 1\right) \geq k\mid
N_{t}\left( 1\right) =i\right) =\frac{\left( \phi _{k}\right) ^{\left(
1\right) }\left( \mathbf{P}\left( N_{t-k}\left( 1\right) =0\right) \right) 
\mathbf{P}\left( N_{t-k}\left( 1\right) =i\right) }{\mathbf{P}\left(
N_{t}\left( 1\right) =i\right) }.  \label{totalsample}
\end{equation}
is the conditional probability distribution of the TMRCA given the whole
population alive at $t$ with size $i$\ is being sampled.
\end{corollary}

Eq. (\ref{totalsample}) has to do with a problem raised in \cite{A} with
some concern on the coalescence time for the whole population of the $t$-th
generation conditioned on the event that it is not extinct. However, Eq. (%
\ref{totalsample}) is somehow more specific as it deals with the coalescence
time for the whole population of the $t$-th generation conditioned on the
event that the population size is precisely $i\geq 2$. \newline

From the preceding result (\ref{totalsample}), upon summing over all
possible values of the current population size, we finally get:\newline

\begin{theorem}
\label{theo2} The distribution of the coalescence time for the whole
population of the $t$-th generation conditioned on the event that it is not
extinct is given by: 
\begin{equation}
\begin{array}{l}
\mathbf{P}\left( \tau _{N_{t}\left( 1\right) ,1}^{\left( t\right) }\left(
1\right) \geq k\mid N_{t}\left( 1\right) >0\right)  \\ 
=\frac{\mathbf{P}\left( N_{t-k}\left( 1\right) >0\right) }{\mathbf{P}\left(
N_{t}\left( 1\right) >0\right) }\left( \phi _{k}\right) ^{\left( 1\right)
}\left( \mathbf{P}\left( N_{t-k}\left( 1\right) =0\right) \right) \mathbf{.}
\end{array}
\label{totalsample2}
\end{equation}
\end{theorem}

\section{\textbf{Proofs}}

From (\ref{1a}), an evaluation of $\mathbf{E}\left( \prod_{l=1}^{j}\frac{%
\left( M_{l}\right) _{i_{l}}}{\left( \sum_{l^{\prime }=1}^{n}M_{l^{\prime
}}\right) _{i_{l}}}\right) $ is needed. In this respect, we have:\newline

\begin{lemma}
$\left( i\right) $\ Let $\left( M_{1},...,M_{n}\right) $\ be $n$\ iid
integral-valued random variables with common pgf $\phi \left( z\right) $.
Then, with $\phi ^{\left( i\right) }\left( z\right) $\ the $i$-th derivative
of $\phi \left( z\right) $\ with respect to $z$, 
\begin{equation*}
\mathbf{E}\left( \frac{\left( M_{1}\right) _{i}}{\left( \sum_{l^{\prime
}=1}^{n}M_{l^{\prime }}\right) _{i}}\right) =\frac{1}{\Gamma \left( i\right) 
}\int_{0}^{1}dz\left( 1-z\right) ^{i-1}\phi \left( z\right) ^{n-1}\phi
^{\left( i\right) }\left( z\right) .
\end{equation*}
$\left( ii\right) $\ With $u\in \left[ 0,1\right] $\ marking the total sum $%
\sum_{l^{\prime }=1}^{n}M_{l^{\prime }},$\ it also holds 
\begin{equation*}
\mathbf{E}\left( \frac{\left( M_{1}\right) _{i}}{\left( \sum_{l^{\prime
}=1}^{n}M_{l^{\prime }}\right) _{i}}u^{\sum_{l^{\prime }=1}^{n}M_{l^{\prime
}}}\right) =\frac{u^{i}}{\Gamma \left( i\right) }\int_{0}^{1}dz\left(
1-z\right) ^{i-1}\phi \left( uz\right) ^{n-1}\phi ^{\left( i\right) }\left(
uz\right) .
\end{equation*}
\end{lemma}

\textbf{Proof}:\newline

$\left( i\right) $ With $\mathcal{M}:=\sum_{l^{\prime }=2}^{n}M_{l^{\prime
}} $, we have 
\begin{equation*}
\mathbf{E}\frac{\left( M_{1}\right) _{i}}{\left( \sum_{l^{\prime
}=1}^{n}M_{l^{\prime }}\right) _{i}}=\sum_{m_{1},m}\mathbf{P}\left(
M_{1}=m_{1}\right) \mathbf{P}\left( \mathcal{M}=m\right) \frac{\left(
m_{1}\right) _{i}}{\left( m_{1}+m\right) _{i}}.
\end{equation*}
On the other hand, $\mathbf{P}\left( \mathcal{M}=m\right) =\left[
z^{m}\right] \phi \left( z\right) ^{n-1}$ and 
\begin{equation*}
\phi ^{\left( i\right) }\left( z\right) =\sum_{m_{1}\geq i}\left(
m_{1}\right) _{i}\mathbf{P}\left( M_{1}=m_{1}\right) z^{m_{1}-i}.
\end{equation*}
Thus, by the beta function identity 
\begin{eqnarray*}
&&\frac{1}{\Gamma \left( i\right) }\int_{0}^{1}dz\left( 1-z\right)
^{i-1}\phi \left( z\right) ^{n-1}\phi ^{\left( i\right) }\left( z\right) \\
&=&\sum_{m_{1},m}\mathbf{P}\left( M_{1}=m_{1}\right) \mathbf{P}\left( 
\mathcal{M}=m\right) \frac{\left( m_{1}\right) _{i}}{\Gamma \left( i\right) }%
\int_{0}^{1}dz\left( 1-z\right) ^{i-1}z^{m_{1}+m-i} \\
&=&\sum_{m_{1},m}\mathbf{P}\left( M_{1}=m_{1}\right) \mathbf{P}\left( 
\mathcal{M}=m\right) \left( m_{1}\right) _{i}\frac{\Gamma \left(
m_{1}+m-i+1\right) }{\Gamma \left( m_{1}+m+1\right) },
\end{eqnarray*}
with $\frac{\Gamma \left( m_{1}+m-i+1\right) }{\Gamma \left(
m_{1}+m+1\right) }=\frac{1}{\left( m_{1}+m\right) _{i}}$. Concerning $\left(
ii\right) $, 
\begin{eqnarray*}
&&\frac{u^{i}}{\Gamma \left( i\right) }\int_{0}^{1}dz\left( 1-z\right)
^{i-1}\phi \left( uz\right) ^{n-1}\phi ^{\left( i\right) }\left( uz\right) \\
&=&\sum_{m_{1},m}\mathbf{P}\left( M_{1}=m_{1}\right) \mathbf{P}\left( 
\mathcal{M}=m\right) \frac{\left( m_{1}\right) _{i}}{\Gamma \left( i\right) }%
u^{m_{1}+m}\int_{0}^{1}dz\left( 1-z\right) ^{i-1}z^{m_{1}+m-i} \\
&=&\sum_{m_{1},m}\mathbf{P}\left( M_{1}=m_{1}\right) \mathbf{P}\left( 
\mathcal{M}=m\right) \frac{\left( m_{1}\right) _{i}}{\left( m_{1}+m\right)
_{i}}u^{m_{1}+m}.\text{ }\Box
\end{eqnarray*}

A slightly extended version of the previous Lemma is:\newline

\begin{corollary}
With $i_{1}+...+i_{j}=i$, it holds 
\begin{equation}
\mathbf{E}\left( \prod_{l=1}^{j}\frac{\left( M_{l}\right) _{i_{l}}}{\left(
\sum_{l^{\prime }=1}^{n}M_{l^{\prime }}\right) _{i_{l}}}\right) =\frac{1}{%
\Gamma \left( i\right) }\int_{0}^{1}dz\left( 1-z\right) ^{i-1}\phi \left(
z\right) ^{n-j}\prod_{l=1}^{j}\phi ^{\left( i_{l}\right) }\left( z\right) .
\label{1b}
\end{equation}

More generally, with $u\in \left[ 0,1\right] $, we have 
\begin{equation}
\begin{array}{l}
\mathbf{E}\left( \prod_{l=1}^{j}\frac{\left( M_{l}\right) _{i_{l}}}{\left(
\sum_{l^{\prime }=1}^{n}M_{l^{\prime }}\right) _{i_{l}}}u^{\sum_{l^{\prime
}=1}^{n}M_{l^{\prime }}}\right)  \\ 
=\frac{u^{i}}{\Gamma \left( i\right) }\int_{0}^{1}dz\left( 1-z\right)
^{i-1}\phi \left( uz\right) ^{n-j}\prod_{l=1}^{j}\phi ^{\left( i_{l}\right)
}\left( uz\right) .
\end{array}
\label{1b2}
\end{equation}
\newline
\end{corollary}

\emph{Proof of Proposition \ref{prop1}.}

Plugging the latter fundamental identity (\ref{1b}) into Eq. (\ref{1a}), we
therefore obtain 
\begin{eqnarray*}
&&\widehat{P}_{i,j}^{\left( t\right) }\left( n_{0}\right) \\
&=&\mathbf{P}\left( \widehat{N}_{t-1}\left( n_{0}\right) =j\mid \widehat{N}%
_{t}\left( n_{0}\right) =i\right) \\
&=&\frac{1}{\Gamma \left( i\right) }\sum_{n\geq j}\binom{n}{j}\mathbf{P}%
\left( N_{t-1}\left( n_{0}\right) =n\right) \int_{0}^{1}dz\left( 1-z\right)
^{i-1}\phi \left( z\right) ^{n-j}\sum_{i_{1}+...+i_{j}=i}^{*}\binom{i}{%
i_{1}...i_{j}}\prod_{l=1}^{j}\phi ^{\left( i_{l}\right) }\left( z\right)
\end{eqnarray*}
and Proposition $1$ follows because formula (\ref{1c}) involves the $j-$th
derivative of \emph{\ }$\phi _{t-1}\left( z\right) ^{n_{0}}:=\sum_{n}z^{n}%
\mathbf{P}\left( N_{t-1}\left( n_{0}\right) =n\right) $, namely $\left[ \phi
_{t-1}^{n_{0}}\right] ^{\left( j\right) }\left( z\right) =\sum_{n\geq
j}\left( n\right) _{j}z^{n-j}\mathbf{P}\left( N_{t-1}\left( n_{0}\right)
=n\right) $, evaluated at $\phi \left( z\right) $. $\Box $\newline

\emph{Proof of Corollary \ref{corol1}.}

The identity (\ref{1c2}) follows from (\ref{1c}) and the fact (discussed in
Remark $1$) that $\mathbf{P}\left( N_{t}\left( n_{0}\right) \geq i\right)
=\sum_{j=1}^{i}\widehat{P}_{i,j}^{\left( t\right) }\left( n_{0}\right) $.
Note that (\ref{1c2}) constitutes of an alternative representation to $%
\mathbf{P}\left( N_{t}\left( n_{0}\right) \geq i\right) =\sum_{j\geq
i}\left[ z^{j}\right] \phi _{t}\left( z\right) ^{n_{0}}.$ $\Box $\newline

\emph{Proof of Theorem \ref{theo1}. }Because the underlying branching
process $\left\{ N_{t}\left( n_{0}\right) \right\} $ is time-homogeneous,
with $N_{k,t}^{\left( l^{\prime }\right) }$ the offspring at time $t$ of the 
$l^{\prime }-$th individual among those in number $N_{k}\left( n_{0}\right) $
alive at time $k$ (so with $\sum_{l^{\prime }=1}^{N_{k}\left( n_{0}\right)
}N_{k,t}^{\left( l^{\prime }\right) }=N_{t}\left( n_{0}\right) $) and
considering the transition probability of $\left\{ \widehat{N}_{s}\left(
n_{0}\right) \right\} $ between steps $t$ and $k<t$ (backwards), we get

\begin{equation*}
\begin{array}{l}
\mathbf{P}\left( \infty >\tau _{i,j}^{\left( t\right) }\left( n_{0}\right)
\geq k\right) \\ 
=\mathbf{E}\left( \binom{N_{k}\left( n_{0}\right) }{j}%
\sum_{i_{1}+...+i_{j}=i}^{*}\binom{i}{i_{1}...i_{j}}\prod_{l=1}^{j}\frac{%
\left( N_{k,t}^{\left( l\right) }\right) _{i_{l}}}{\left( \sum_{l^{\prime
}=1}^{N_{k}\left( n_{0}\right) }N_{k,t}^{\left( l^{\prime }\right) }\right)
_{i_{l}}}\right) \\ 
=\frac{1}{j!\Gamma \left( i\right) }\sum_{i_{1}+...+i_{j}=i}^{*}\binom{i}{%
i_{1}...i_{j}}\int_{0}^{1}dz\left( 1-z\right) ^{i-1}\left[ \phi
_{k}^{n_{0}}\right] ^{\left( j\right) }\left( \phi _{t-k}\left( z\right)
\right) \prod_{l=1}^{j}\phi _{t-k}^{\left( i_{l}\right) }\left( z\right) .
\end{array}
\end{equation*}

The last equation is obtained from (\ref{1a}), (\ref{1b}) and (\ref{1c})
while substituting $\phi _{t-k}\left( z\right) $ to $\phi \left( z\right) $
and $\phi _{k}\left( z\right) $ to $\phi _{t-1}\left( z\right) .$ $\Box $%
\newline

\emph{Proof of Corollary \ref{corol2}.}

When $j=1,$%
\begin{equation*}
\begin{array}{l}
\mathbf{P}\left( \infty >\tau _{i,1}^{\left( t\right) }\left( n_{0}\right)
\geq k\right) =\mathbf{E}\left( N_{k}\left( n_{0}\right) \frac{\left(
N_{k,t}\right) _{i}}{\left( \sum_{l^{\prime }=1}^{N_{k}\left( n_{0}\right)
}N_{k,t}^{\left( l^{\prime }\right) }\right) _{i}}\right) \\ 
=\frac{1}{\Gamma \left( i\right) }\int_{0}^{1}dz\left( 1-z\right)
^{i-1}\left( \phi _{k}^{n_{0}}\right) ^{\left( 1\right) }\left( \phi
_{t-k}\left( z\right) \right) \phi _{t-k}^{\left( i\right) }\left( z\right)
\\ 
=\frac{1}{\Gamma \left( i\right) }\int_{0}^{1}dz\left( 1-z\right)
^{i-1}\left[ \phi _{t}^{n_{0}}\left( z\right) \right] ^{\left( 1\right) }%
\frac{\phi _{t-k}^{\left( i\right) }\left( z\right) }{\phi _{t-k}^{\left(
1\right) }\left( z\right) }.
\end{array}
\end{equation*}
The later equality makes use of the observation that $\left[ \phi _{t}\left(
z\right) ^{n_{0}}\right] ^{\left( 1\right) }=\left[ \phi _{k}\left( \phi
_{t-k}\left( z\right) \right) ^{n_{0}}\right] ^{\left( 1\right) }=\left[
\phi _{k}^{n_{0}}\right] ^{\left( 1\right) }\left( \phi _{t-k}\left(
z\right) \right) \phi _{t-k}^{\left( 1\right) }\left( z\right) $. The final
proof of (\ref{1e}) and (\ref{1f}) follows from (\ref{1b2}) and (\ref{1d})
with $j=1$ \footnote{%
When $i=2$, we recover Corollary $1$ in \cite{Lam}, which was derived
differently.}. $\Box $\newline

\emph{Proof of Proposition \ref{prop2}.}

(\ref{finite}) is obtained while putting $k=0$ in (\ref{1de}).

(\ref{ip}) is obtained upon integrating (\ref{finite}) by parts.\newline

\emph{Proof of Corollary \ref{corol3}.}

The identity (\ref{I2}) follows from plugging $i=2$ in the integration by
parts formula (\ref{ip}) and observing $\mathbf{P}\left( \tau _{1,1}^{\left(
t\right) }\left( n_{0}\right) <\infty \right) =1-\mathbf{P}\left(
N_{t}\left( n_{0}\right) =0\right) $. For this point, see also Corollary $1$
in \cite{Lam}. $\Box $\newline

\emph{Proof of Corollary \ref{corol4}.}

The identity (\ref{totalsample}) follows from (\ref{1f}) plugging $j=i$ and $%
n_{0}=1$. Furthermore, 
\begin{equation*}
\mathbf{P}\left( \tau _{i,1}^{\left( t\right) }\left( 1\right) \geq k\mid
N_{t}\left( 1\right) =i\right) =\frac{\left[ u^{0}\right]
\int_{0}^{1}dz\left( 1-z\right) ^{i-1}\left( \phi _{k}\right) ^{\left(
1\right) }\left( \phi _{t-k}\left( uz\right) \right) \phi _{t-k}^{\left(
i\right) }\left( uz\right) }{\Gamma \left( i\right) \mathbf{P}\left(
N_{t}\left( 1\right) =i\right) }
\end{equation*}
\begin{equation*}
=\frac{\left( \phi _{k}\right) ^{\left( 1\right) }\left( \phi _{t-k}\left(
0\right) \right) \phi _{t-k}^{\left( i\right) }\left( 0\right) }{\Gamma
\left( i+1\right) \mathbf{P}\left( N_{t}\left( 1\right) =i\right) }=\frac{%
\left( \phi _{k}\right) ^{\left( 1\right) }\left( \mathbf{P}\left(
N_{t-k}\left( 1\right) =0\right) \right) \mathbf{P}\left( N_{t-k}\left(
1\right) =i\right) }{\mathbf{P}\left( N_{t}\left( 1\right) =i\right) }.\text{
}\Box
\end{equation*}
\newline
\emph{Remark:} Because Eq. (\ref{totalsample}) is not straightforward, let
us illustrate it in a very simple case. Suppose $\phi \left( z\right) =\pi
_{0}+\pi _{1}z+\pi _{2}z^{2}$ (a random binary tree). Suppose $t=2$ and $i=2$
in (\ref{totalsample}). The event $N_{2}\left( 1\right) =2$ may occur
because the ancestor has $2$ offsprings, each generating next a single
offspring and then $\tau _{2,1}^{\left( 2\right) }\left( 1\right) =0$ with
probability $\pi _{2}\pi _{1}^{2}/\mathbf{P}\left( N_{2}\left( 1\right)
=2\right) $. It may also occur because the ancestor has $1$ offspring,
itself generating $2$ offspring (an event with probability $\pi _{1}\pi _{2}$%
) or because the ancestor has $2$ offspring, only one of which generating
two offspring (an event with probability $2\pi _{2}\pi _{0}\pi _{2}$) and
then one expects $\tau _{2,1}^{\left( 2\right) }\left( 1\right) =1$ with
probability $\left( \pi _{1}\pi _{2}+2\pi _{2}\pi _{0}\pi _{2}\right) /%
\mathbf{P}\left( N_{2}\left( 1\right) =2\right) $ with $\mathbf{P}\left(
N_{2}\left( 1\right) =2\right) =\left[ z^{2}\right] \phi _{2}\left( z\right)
=\pi _{2}\pi _{1}^{2}+$ $\pi _{1}\pi _{2}+2\pi _{2}\pi _{0}\pi _{2}$.
Putting $t=2$, $i=2$ and $k=1$ in (\ref{totalsample}), with $N_{1}\left(
1\right) =M$ and $\left( \phi _{1}\right) ^{\left( 1\right) }\left( z\right)
=\pi _{1}+2\pi _{2}z,$%
\begin{eqnarray*}
\mathbf{P}\left( \tau _{2,1}^{\left( 2\right) }\left( 1\right) =1\mid
N_{2}\left( 1\right) =2\right) &=&\frac{\left( \phi _{1}\right) ^{\left(
1\right) }\left( \mathbf{P}\left( N_{1}\left( 1\right) =0\right) \right) 
\mathbf{P}\left( N_{1}\left( 1\right) =2\right) }{\mathbf{P}\left(
N_{2}\left( 1\right) =2\right) } \\
&=&\frac{\left( \pi _{1}+2\pi _{2}\pi _{0}\right) \pi _{2}}{\mathbf{P}\left(
N_{2}\left( 1\right) =2\right) },
\end{eqnarray*}
which is the expected result. And, because $\mathbf{P}\left( N_{2}\left(
1\right) =2\right) =\pi _{2}\pi _{1}^{2}+$ $\pi _{1}\pi _{2}+2\pi _{2}\pi
_{0}\pi _{2},$%
\begin{eqnarray*}
&&\mathbf{P}\left( \tau _{2,1}^{\left( 2\right) }\left( 1\right) =0\mid
N_{2}\left( 1\right) =2\right) \\
&=&\mathbf{P}\left( \tau _{2,1}^{\left( 2\right) }\left( 1\right) \geq 0\mid
N_{2}\left( 1\right) =2\right) -\mathbf{P}\left( \tau _{2,1}^{\left(
2\right) }\left( 1\right) =1\mid N_{2}\left( 1\right) =2\right) \\
&=&1-\mathbf{P}\left( \tau _{2,1}^{\left( 2\right) }\left( 1\right) =1\mid
N_{2}\left( 1\right) =2\right) =\frac{\pi _{2}\pi _{1}^{2}}{\mathbf{P}\left(
N_{2}\left( 1\right) =2\right) },
\end{eqnarray*}
again the expected result.\newline

\emph{Proof of Theorem \ref{theo2}. }From the result (\ref{totalsample}),
upon summing over all possible values of the current population size, we
get: 
\begin{equation*}
\begin{array}{l}
\mathbf{P}\left( \tau _{N_{t}\left( 1\right) ,1}^{\left( t\right) }\left(
1\right) \geq k\mid N_{t}\left( 1\right) >0\right) \\ 
=\frac{1}{\mathbf{P}\left( N_{t}\left( 1\right) >0\right) }\sum_{i\geq 1}%
\mathbf{P}\left( \tau _{i,1}^{\left( t\right) }\left( 1\right) \geq k\mid
N_{t}\left( 1\right) =i\right) \mathbf{P}\left( N_{t}\left( 1\right)
=i\right) \\ 
=\frac{\mathbf{P}\left( N_{t-k}\left( 1\right) >0\right) }{\mathbf{P}\left(
N_{t}\left( 1\right) >0\right) }\left( \phi _{k}\right) ^{\left( 1\right)
}\left( \mathbf{P}\left( N_{t-k}\left( 1\right) =0\right) \right) \mathbf{.}
\end{array}
\end{equation*}

\section{\textbf{Special reproduction laws examples}}

\subsection{The b-ary tree}

Before proceeding with more meaningful reproduction laws, it is of interest
to first consider the deterministic $b$-ary tree case for which $\phi \left(
z\right) =z^{b}$ where $b\geq 2$ is an integer. This is a trivial example of
a strictly supercritical branching process, exploding with probability $1$
in infinite time. In this case, $\phi _{t}\left( z\right) =z^{b^{t}}$ and $%
N_{t}\left( n_{0}\right) =n_{0}b^{t}\geq 2$ as $t\geq 0.$

\subsubsection{\textbf{Coalescence times for }$i$\textbf{\ sampled
individuals}}

In this section, we use the formulas derived in the first half of this paper
to derive the probability of $i$ individuals to have a common ancestor, some
transition matrix elements (the $\widehat{P}_{i,1}^{\left( t\right) }\left(
n_{0}\right) $) and the conditional probability on the TMRCA given that it
is finite: 
\begin{equation*}
\mathbf{P}\left( \tau _{i,1}^{\left( t\right) }\left( n_{0}\right) <\infty
\right) =\frac{n_{0}\left( b^{t}\right) _{i}}{\Gamma \left( i\right) }%
\int_{0}^{1}dz\left( 1-z\right) ^{i-1}z^{\left( n_{0}-1\right)
b^{t}+b^{t}-i}=\frac{n_{0}\left( b^{t}\right) _{i}}{\left( n_{0}b^{t}\right)
_{i}}.
\end{equation*}

Note that, as required, $\mathbf{P}\left( \tau _{i,1}^{\left( t\right)
}\left( n_{0}\right) =\infty \right) =1$ if $i>b^{t}$. As $n_{0}\rightarrow
\infty $, $\mathbf{P}\left( \tau _{i,1}^{\left( t\right) }\left(
n_{0}\right) <\infty \right) \sim n_{0}^{-\left( i-1\right) }\left(
b^{t}\right) _{i}b^{-it}$. 
\begin{equation*}
\widehat{P}_{i,1}^{\left( t\right) }\left( n_{0}\right) =\mathbf{P}\left(
\infty >\tau _{i,1}^{\left( t\right) }\left( n_{0}\right) \geq t-1\right) =%
\frac{n_{0}b^{t-1}\left( b\right) _{i}}{\Gamma \left( i\right) \left(
n_{0}b^{t}\right) _{i}}\text{ and}
\end{equation*}
\begin{equation*}
\mathbf{P}\left( \tau _{i,1}^{\left( t\right) }\left( n_{0}\right) \geq
k\mid \tau _{i,1}^{\left( t\right) }\left( n_{0}\right) <\infty \right) =%
\frac{b^{k}\left( b^{t-k}\right) _{i}}{\Gamma \left( i\right) \left(
b^{t}\right) _{i}},
\end{equation*}
independently of $n_{0}$. Note that the one-step ratio of a triple to a
binary merger obeys 
\begin{equation*}
\frac{\widehat{P}_{3,1}^{\left( t\right) }\left( n_{0}\right) }{\widehat{P}%
_{2,1}^{\left( t\right) }\left( n_{0}\right) }=\frac{b-2}{n_{0}b^{t}-2}\sim 
\frac{b-2}{n_{0}b^{t}}\rightarrow 0\text{ as }n_{0}\rightarrow \infty \text{
or }t\rightarrow \infty .\text{ }
\end{equation*}
The three first results derive from Proposition \ref{prop2}, Corollary \ref
{corol2} and Corollary \ref{corol4}.

\subsubsection{\textbf{Coalescence time for pairs of tips}}

When\textbf{\ }$i=2$, 
\begin{equation*}
\mathbf{P}\left( \infty >\tau _{2,1}^{\left( t\right) }\left( n_{0}\right)
\geq k\right) =\frac{b^{t-k}-1}{n_{0}b^{t}-1}.
\end{equation*}
In this non-random case, we get 
\begin{eqnarray*}
\mathbf{P}\left( N_{t}\left( n_{0}\right) \geq 2,\tau _{2,1}^{\left(
t\right) }\left( n_{0}\right) =\infty \right) &=&\mathbf{P}\left( \tau
_{2,1}^{\left( t\right) }\left( n_{0}\right) =\infty \right) =\frac{\left(
n_{0}-1\right) b^{t}}{n_{0}b^{t}-1}, \\
\mathbf{P}\left( \tau _{2,1}^{\left( t\right) }\left( n_{0}\right) <\infty
\right) &=&\frac{b^{t}-1}{n_{0}b^{t}-1}\underset{t\rightarrow \infty }{%
\rightarrow }n_{0}^{-1}.
\end{eqnarray*}
When $n_{0}$ is large, the probability that a sampled random pair at $t$ has
a common ancestor (is a member of the family of the same founder) is of
order $n_{0}^{-1}\left( 1-b^{-t}\right) .$ If in addition $t$ is large, this
probability is of order $n_{0}^{-1}.$ Furthermore, 
\begin{equation*}
\mathbf{P}\left( \tau _{2,1}^{\left( t\right) }\left( n_{0}\right) =k\mid
\tau _{2,1}^{\left( t\right) }\left( n_{0}\right) <\infty \right) =\frac{%
b^{-k}\left( 1-b^{-1}\right) }{1-b^{-t}}\text{, }k=0,...,t-1,
\end{equation*}
which is a truncated geometric$\left( b^{-1}\right) $ distribution
approaching the geometric distribution when $t$ itself gets large.

The distribution of the coalescence time for the whole population of the $t$%
-th generation in this deterministic case is uninteresting, because the law
of $\tau _{N_{t}\left( n_{0}\right) ,1}^{\left( t\right) }\left(
n_{0}\right) $ clearly is a Dirac mass at generation $0$.

As in this simplistic example, whenever $\phi _{t}\left( z\right) $ is
explicitly known, $\mathbf{P}\left( \tau _{2,1}^{\left( t\right) }\left(
n_{0}\right) \geq k\right) $ for instance can be explicitly computed in
principle. Let us discuss two cases: the first one corresponding to $\mu =%
\mathbf{E}\left( M\right) <\infty $, the other one to $\mu =\infty $. The
first example with finite mean illustrates what one may expect when the
branching process is either (sub-)critical or supercritical. The second one
deals with an example with infinite mean.

\subsection{Linear-fractional mechanism with finite mean $\mu $}

Let $p_{0},p$ be probabilities and put $q_{0}=1-p_{0},$ $q=1-p.$ The
branching mechanism 
\begin{equation*}
\phi \left( z\right) =q_{0}+p_{0}qz/\left( 1-pz\right) =q_{0}+p_{0}z/\left(
1+p/q\left( 1-z\right) \right)
\end{equation*}
is called linear fractional with $\mu =p_{0}/q$. Introduce $a,b>0$ defined
by $p_{0}=1/\left( a+b\right) $ and $p=b/\left( a+b\right) $ so that $\mu
=1/a$. The supercritical case $\mu >1$ (respectively subcritical case $\mu
<1 $) corresponds to $a<1$ (respectively $a>1$). Then 
\begin{eqnarray*}
\phi \left( z\right) &=&\left( 1-\frac{1}{a+b}\right) +\frac{1}{a+b}\frac{z}{%
1+\frac{b}{a}\left( 1-z\right) }\text{ and} \\
\phi _{t}\left( z\right) &=&\left( 1-\frac{1}{a_{t}+b_{t}}\right) +\frac{1}{%
a_{t}+b_{t}}\frac{z}{1+\frac{b_{t}}{a_{t}}\left( 1-z\right) },
\end{eqnarray*}
where $a_{t}=a^{t}$ and $b_{t}=b\left( 1+a+...+a^{t-1}\right) $, $\mu
_{t}=a^{-t}.$

Note $a+b=1\Rightarrow a_{t}+b_{t}=1$ and $\phi \left( 0\right) =\mathbf{P}%
\left( M=0\right) =\phi _{t}\left( 0\right) =\mathbf{P}\left( N_{t}\left(
1\right) =0\right) =\mathbf{P}\left( N_{t}\left( n_{0}\right) =0\right) =0$.
Clearly, 
\begin{equation*}
\phi _{t}\left( z\right) =\frac{A_{t}-B_{t}z}{C_{t}-D_{t}z},
\end{equation*}
where 
\begin{equation*}
A_{t}=a_{t}+b_{t}-1,\text{ }B_{t}=b_{t}-1,\text{ }C_{t}=a_{t}+b_{t},\text{ }%
D_{t}=b_{t}.
\end{equation*}
In the sequel, to avoid additional complexities, we shall limit ourselves to
the geometrical case, which is a special case of the linear-fractional
model. The linear-fractional case can be treated along similar lines.

\subsubsection{\textbf{The geometrical model}}

If $p_{0}=p$, then $\phi \left( z\right) =q+pqz/\left( 1-pz\right) =q/\left(
1-pz\right) $, the classical geometric distribution with mean $\mu =p/q$.
Here, $\phi \left( 0\right) =\mathbf{P}\left( M=0\right) =q\neq 0$ and $%
a=q/p=1/\mu $, $b=1$. Furthermore,

- If $\mu \neq 1$: $\phi _{t}\left( 0\right) =\mathbf{P}\left( N_{t}\left(
1\right) =0\right) =1-\frac{1}{a_{t}+b_{t}}$, $a_{t}=a^{t}$ and $%
b_{t}=1+a+...+a^{t-1}=\left( a^{t}-1\right) /\left( a-1\right) $.

- If $\mu =1$: $\phi _{t}\left( 0\right) =\mathbf{P}\left( N_{t}\left(
1\right) =0\right) =1-\frac{1}{1+t}$, $a_{t}=1$ and $b_{t}=t$.

Let us illustrate the peculiarities of the geometric model. In this case,
with $\phi \left( z\right) =q/\left( 1-pz\right) $,

- If $\mu \neq 1$: 
\begin{equation*}
\phi _{t}\left( z\right) =\frac{A_{t}-B_{t}z}{C_{t}-D_{t}z}=\alpha _{t}\frac{%
1-\beta _{t}z}{1-\delta _{t}z},
\end{equation*}
where $a=q/p=\mu ^{-1}$, $a_{t}=a^{t}$, $b_{t}=1+a+...+a^{t-1}=\left(
a^{t}-1\right) /\left( a-1\right) $, 
\begin{equation*}
A_{t}=a_{t}+b_{t}-1,\text{ }B_{t}=b_{t}-1,\text{ }C_{t}=a_{t}+b_{t},\text{ }%
D_{t}=b_{t}\text{ and}
\end{equation*}
\begin{equation*}
\alpha _{t}=\frac{A_{t}}{C_{t}}=a\frac{a^{t}-1}{a^{t+1}-1};\text{ }\beta
_{t}=\frac{B_{t}}{A_{t}}=\frac{a^{t-1}-1}{a^{t}-1}\text{ and }\delta _{t}=%
\frac{D_{t}}{C_{t}}=\frac{a^{t}-1}{a^{t+1}-1}=\beta _{t+1}.
\end{equation*}

- If $\mu =1$: $a=a^{t}=1$, $b_{t}=t$, $A_{t}=t,$ $B_{t}=t-1,$ $C_{t}=1+t,$ $%
D_{t}=t$, 
\begin{equation*}
\alpha _{t}=\frac{A_{t}}{C_{t}}=\frac{t}{1+t};\text{ }\beta _{t}=\frac{B_{t}%
}{A_{t}}=\frac{t-1}{t}\text{ and }\delta _{t}=\frac{D_{t}}{C_{t}}=\frac{t}{%
1+t}\text{ and}
\end{equation*}
\begin{equation*}
\phi _{t}\left( z\right) =\frac{A_{t}-B_{t}z}{C_{t}-D_{t}z}=\alpha _{t}\frac{%
1-\beta _{t}z}{1-\delta _{t}z}.
\end{equation*}
In any case, we have 
\begin{equation*}
\phi _{t}^{\prime }\left( z\right) =\alpha _{t}\frac{\delta _{t}-\beta _{t}}{%
\left( 1-\delta _{t}z\right) ^{2}}\text{ and }\phi _{t}^{^{\prime \prime
}}\left( z\right) =2\alpha _{t}\delta _{t}\left( \delta _{t}-\beta
_{t}\right) \left( 1-\delta _{t}z\right) ^{-3}.
\end{equation*}
With these simple preliminaries at hand, let us investigate some coalescence
times results pertaining to the geometric branching mechanism.

\subsubsection{\textbf{Tail probability of }$\tau _{2,1}^{\left( t\right)
}\left( n_{0}\right) $\textbf{\ and Lauricella functions}}

We express the tail probability of the TMRCA using a Lauricella
hypergeometric function. We have 
\begin{eqnarray*}
&&\mathbf{P}\left( \infty >\tau _{2,1}^{\left( t\right) }\left( n_{0}\right)
\geq k\right) \\
&=&n_{0}\delta _{t-k}\left( \delta _{t}-\beta _{t}\right) \alpha
_{t}^{n_{0}}\times L^{\left( 3\right) }\left( 1;1-n_{0},1,n_{0}+1;3;\beta
_{t},\delta _{t-k},\delta _{t}\right) ,
\end{eqnarray*}
where $L^{\left( 3\right) }$ is a fourth-kind Lauricella hypergeometric
function of order $3$, \cite{Carl} and \cite{Lau}. That result follows from
Corollary \ref{corol2}. $L^{\left( 3\right) }$ is defined by 
\begin{eqnarray*}
L^{\left( 3\right) }\left( 1;1-n_{0},1,n_{0}+1;3;\beta _{t},\delta
_{t-k},\delta _{t}\right) &=&\sum_{m\geq 0}\frac{\left[ 1\right] _{m}}{%
\left[ 3\right] _{m}}\left[ z^{m}\right] \frac{1}{\left( 1-\delta
_{t-k}z\right) }\frac{\left( 1-\beta _{t}z\right) ^{n_{0}-1}}{\left(
1-\delta _{t}z\right) ^{n_{0}+1}} \\
&=&\sum_{m\geq 0}\frac{2}{\left( m+1\right) \left( m+2\right) }\Lambda _{m},
\end{eqnarray*}
where 
\begin{equation*}
\Lambda _{m}=\left[ z^{m}\right] \frac{1}{\left( 1-\delta _{t-k}z\right) }%
\frac{\left( 1-\beta _{t}z\right) ^{n_{0}-1}}{\left( 1-\delta _{t}z\right)
^{n_{0}+1}}\text{.}
\end{equation*}

Putting $k=0$ in (\ref{closfm}), we get 
\begin{eqnarray*}
\mathbf{P}\left( \tau _{2,1}^{\left( t\right) }\left( n_{0}\right) <\infty
\right) &=&2n_{0}\delta _{t}\left( \delta _{t}-\beta _{t}\right) \alpha
_{t}^{n_{0}}\int_{0}^{1}dz\left( 1-z\right) \frac{\left( 1-\beta
_{t}z\right) ^{n_{0}-1}}{\left( 1-\delta _{t}z\right) ^{n_{0}+2}} \\
&=&n_{0}\delta _{t}\left( \delta _{t}-\beta _{t}\right) \alpha
_{t}^{n_{0}}\cdot F_{1}\left( 1;1-n_{0},n_{0}+2;3;\beta _{t},\delta
_{t}\right) ,
\end{eqnarray*}
where $F_{1}$ is an order-$2$ Lauricella function, also called an Appell
hypergeometric function. \newline

In the sequel, we shall need the following trivial observations: when $t$
gets large and in the subcritical case $a>1$ ($\mu <1$), $\delta _{t}\sim
\mu \left( 1-\mu ^{t}\left( 1-\mu \right) \right) $ and $\beta _{t}\sim \mu
\left( 1-\mu ^{t-1}\left( 1-\mu \right) \right) $, approaching $\mu $, with $%
\delta _{t}-\beta _{t}\sim \mu ^{t}\left( 1-\mu \right) ^{2}$, whereas in
the supercritical case $a<1$ ($\mu >1$), $\beta _{t}\sim 1-a^{t-1}\left(
1-a\right) $ and $\delta _{t}\sim 1-a^{t}\left( 1-a\right) $, approaching $%
1, $ with $\delta _{t}-\beta _{t}\sim \left( 1-a\right) ^{2}a^{t-1}.$ In the
critical case, $\delta _{t}\sim 1-1/\left( t+1\right) $ and $\beta _{t}\sim
1-1/t,$ with $\delta _{t}-\beta _{t}\sim 1/\left[ t\left( t+1\right) \right]
.$

\subsubsection{\textbf{A single founder }$n_{0}=1$}

Because the computations of the law of $\tau _{2,1}^{\left( t\right) }\left(
n_{0}\right) $ with $n_{0}\neq 1$ are rather cumbersome, we shall completely
perform the computations only in the simpler case $n_{0}=1$. In this case,
using identities provided by $\emph{Mathematica}$, we have the simple
closed-form formula: 
\begin{equation}
\mathbf{P}\left( \infty >\tau _{2,1}^{\left( t\right) }\left( 1\right) \geq
k\right) =\frac{2\alpha _{t}\left( \delta _{t}-\beta _{t}\right) \delta
_{t-k}\left[ \delta _{t-k}-\delta _{t}+\left( 1-\delta _{t-k}\right) \log 
\frac{1-\delta _{t-k}}{1-\delta _{t}}\right] }{\left( \delta _{t}-\delta
_{t-k}\right) ^{2}}.  \label{closfm}
\end{equation}

The formula (\ref{closfm}) allows to evaluate large-$t$ asymptotics in the
three regimes: sub-critical, critical and super-critical. Specifically:

\begin{itemize}
\item  If $\mu <1$ (sub-critical): 
\begin{equation}
\mathbf{P}\left( t-\tau _{2,1}^{\left( t\right) }\left( 1\right) \leq l\mid
\tau _{2,1}^{\left( t\right) }\left( 1\right) <\infty \right) \sim \mu -%
\frac{\mu }{3}\mu ^{l}\left( 3-2\mu \right) ,\text{ }l\in \left\{
1,...,t\right\} ,  \label{e1}
\end{equation}
showing that coalescence for pairs of tips occurs in the near past.
\end{itemize}

This result follows from using the explicit expressions of $\alpha _{t},$ $%
\beta _{t}$ and $\delta _{t}$, 
\begin{eqnarray*}
&&\mathbf{P}\left( \infty >\tau _{2,1}^{\left( t\right) }\left( 1\right)
\geq k\right) \\
&&\underset{\text{large }t}{\sim }\mu \left( 1-\mu \right) \mu ^{t}-\frac{1}{%
3}\left( 1-\mu \right) \mu ^{1-k}\left( 3-2\mu -4\mu ^{k+1}\right) \mu ^{2t}.
\end{eqnarray*}
Similarly, 
\begin{equation*}
\mathbf{P}\left( \tau _{2,1}^{\left( t\right) }\left( 1\right) <\infty
\right) =\frac{a-1}{a^{t+1}-1}\sim \left( 1-\mu \right) \mu ^{t}.
\end{equation*}
This probability decays geometrically as $\mu ^{t}$ when $t$ gets large. The
result follows from these two formulas.

\begin{itemize}
\item  If $\mu >1$ (supercritical): with $k\in \Bbb{N}_{0}$, 
\begin{equation}
\mathbf{P}\left( \tau _{2,1}^{\left( t\right) }\left( 1\right) \geq k\mid
\tau _{2,1}^{\left( t\right) }\left( 1\right) <\infty \right) \rightarrow 
\overline{\pi }_{k}:=\frac{2\mu ^{-k}\left( k\log \mu -\left( 1-\mu
^{-k}\right) \right) }{\left( 1-\mu ^{-k}\right) ^{2}},  \label{e3}
\end{equation}
showing that coalescence for pairs occurs in the remote past in the
supercritical regime. For the limiting tail distribution to the right, one
can check that $\overline{\pi }_{k}\rightarrow 1$ as $k\rightarrow 0.$

This result follows from (\ref{closfm}): 
\begin{equation*}
\mathbf{P}\left( \infty >\tau _{2,1}^{\left( t\right) }\left( 1\right) \geq
k\right) \underset{\text{large }t}{\rightarrow }\frac{2\left( \mu -1\right)
\mu ^{-k-1}\left( k\log \mu +\mu ^{-k}-1\right) }{\left( 1-\mu ^{-k}\right)
^{2}},
\end{equation*}
the convergence being at geometric rate $\mu ^{-\left( t+1\right) }$ when $t$
gets large. Moreover, 
\begin{equation*}
\mathbf{P}\left( \tau _{2,1}^{\left( t\right) }\left( 1\right) <\infty
\right) =\frac{1-a}{1-a^{t+1}}\sim \frac{\mu -1}{\mu }\left( 1-\mu ^{-\left(
t+1\right) }\right) \rightarrow \frac{\mu -1}{\mu }.
\end{equation*}
This probability tends to a limit at geometric rate $\mu ^{-\left(
t+1\right) }$ when $t$ gets large.

\item  If $\mu =1$ (critical): 
\begin{equation}
\mathbf{P}\left( \frac{\tau _{2,1}^{\left( t\right) }\left( 1\right) }{t}%
\geq x\mid \tau _{2,1}^{\left( t\right) }\left( 1\right) <\infty \right)
\rightarrow -\frac{2}{x^{2}}\left( x\left( 1-x\right) +\left( 1-x\right)
\log \left( 1-x\right) \right) .  \label{e4}
\end{equation}
In the critical regime, coalescence for pairs occurs in-between the recent
and the remote past.

This result follows from using the explicit expressions of $\alpha _{t},$ $%
\beta _{t}$ and $\delta _{t}$: 
\begin{equation*}
\mathbf{P}\left( \infty >\tau _{2,1}^{\left( t\right) }\left( 1\right) \geq
k\right) =\frac{2\left( t-k\right) }{k^{2}\left( t+1\right) }\left( \left(
t+1\right) \log \frac{t+1}{t-k+1}-k\right) .
\end{equation*}
When $k=\left[ tx\right] $ for some $x\in \left( 0,1\right) $, as $t$ gets
large, 
\begin{equation*}
\mathbf{P}\left( \frac{\tau _{2,1}^{\left( t\right) }\left( 1\right) }{t}%
\geq x\right) \rightarrow \frac{-2\left( 1-x\right) }{tx^{2}}\left( x+\log
\left( 1-x\right) \right) .
\end{equation*}
Moreover, 
\begin{equation*}
\mathbf{P}\left( \tau _{2,1}^{\left( t\right) }\left( 1\right) <\infty
\right) =\frac{t}{\left( 1+t\right) ^{2}}.
\end{equation*}
This probability decays algebraically as $t^{-1}$ when $t$ gets large.
\end{itemize}

\subsubsection{\textbf{Existence of a common ancestor when }$n_{0}>1$}

Formula (\ref{closfm}) with $k=0$ gives the probability that a common
ancestor exists.

When $n_{0}>1$ is fixed, extracting the $\left[ z^{m}\right] -$coefficient
of $\left( 1-\beta _{t}z\right) ^{n_{0}-1}/\left( 1-\delta _{t}z\right)
^{n_{0}+2}$ in the Lauricella integral, we obtain the identity 
\begin{equation*}
\int_{0}^{1}dz\left( 1-z\right) \frac{\left( 1-\beta _{t}z\right) ^{n_{0}-1}%
}{\left( 1-\delta _{t}z\right) ^{n_{0}+2}}=\frac{\left( \beta _{t}-1\right)
-n_{0}\left( \delta _{t}-\beta _{t}\right) +\frac{\left( 1-\beta _{t}\right)
^{n_{0}+1}}{\left( 1-\delta _{t}\right) ^{n_{0}}}}{n_{0}\left(
n_{0}+1\right) \left( \delta _{t}-\beta _{t}\right) ^{2}}.
\end{equation*}
We thus get:

\begin{itemize}
\item  If $\mu <1$ (subcritical): 
\begin{equation}
\mathbf{P}\left( \tau _{2,1}^{\left( t\right) }\left( n_{0}\right) <\infty
\right) \sim n_{0}\mu \left( 1-\mu \right) \mu ^{t}.  \label{e5}
\end{equation}
This probability decays geometrically like $\mu ^{t}$.\newline

\item  If $\mu >1$ (supercritical): 
\begin{equation}
\mathbf{P}\left( \tau _{2,1}^{\left( t\right) }\left( n_{0}\right) <\infty
\right) =2\delta _{t}\alpha _{t}^{n_{0}}\frac{\left( \beta _{t}-1\right)
-n_{0}\left( \delta _{t}-\beta _{t}\right) +\frac{\left( 1-\beta _{t}\right)
^{n_{0}+1}}{\left( 1-\delta _{t}\right) ^{n_{0}}}}{\left( n_{0}+1\right)
\left( \delta _{t}-\beta _{t}\right) }.  \label{e6}
\end{equation}
To the leading order in $t$, we get 
\begin{equation*}
\mathbf{P}\left( \tau _{2,1}^{\left( t\right) }\left( n_{0}\right) <\infty
\right) \rightarrow 2\frac{\mu \left( 1-\mu ^{-n_{0}}\right) -n_{0}\left(
\mu -1\right) \mu ^{-n_{0}}}{\left( \mu -1\right) \left( n_{0}+1\right) },
\end{equation*}
at geometric rate $\mu ^{-t}.$ This probability decays geometrically like $%
\mu ^{-t}$.\newline

\item  If $\mu =1$ (critical): 
\begin{equation}
\begin{array}{l}
\mathbf{P}\left( \tau _{2,1}^{\left( t\right) }\left( n_{0}\right) <\infty
\right) =\frac{2t}{\left( n_{0}+1\right) \left( t+1\right) ^{n_{0}+1}}\left[
\left( t+1\right) ^{n_{0}+1}-t^{n_{0}}\left( n_{0}+t+1\right) \right] \\ 
\sim n_{0}/t.
\end{array}
\label{e7}
\end{equation}
This probability decays algebraically like $t^{-1}$. Using

\begin{equation*}
\mathbf{P}\left( N_{t}\left( n_{0}\right) \geq 2,\tau _{2,1}^{\left(
t\right) }\left( n_{0}\right) =\infty \right) =\mathbf{P}\left( N_{t}\left(
n_{0}\right) \geq 2\right) -\mathbf{P}\left( \tau _{2,1}^{\left( t\right)
}\left( n_{0}\right) <\infty \right) ,
\end{equation*}
and because 
\begin{equation*}
\mathbf{P}\left( N_{t}\left( 1\right) =0\right) =\frac{A_{t}}{C_{t}}=\alpha
_{t}\text{ and}
\end{equation*}
\begin{equation}
\mathbf{P}\left( N_{t}\left( 1\right) =i\right) =\frac{D_{t}^{i-1}\left(
A_{t}D_{t}-B_{t}C_{t}\right) }{C_{t}^{i+1}}=\alpha _{t}\delta
_{t}^{i-1}\left( \delta _{t}-\beta _{t}\right) ,\text{ }i\geq 1  \label{e8}
\end{equation}
are known, $\mathbf{P}\left( N_{t}\left( n_{0}\right) \geq 2\right)
=1-\alpha _{t}\left( 1+\delta _{t}-\beta _{t}\right) $ and $\mathbf{P}\left(
N_{t}\left( n_{0}\right) \geq 2,\tau _{2,1}^{\left( t\right) }\left(
n_{0}\right) =\infty \right) $ follows. This is the probability that no
coalescence occurs because the 2 sampled individuals are not related to the
same founder. The large-$t$ estimate of these probabilities can easily be
derived in the (sub)-critical, supercritical cases. We skip the details.
\end{itemize}

\subsubsection{\textbf{Coalescence time for the whole population of the }$t$%
\textbf{-th generation conditioned on its size and the event that it is not
extinct}}

We come back to $n_{0}=1$ for simplicity. If $\mu \neq 1$, following
Corollary \ref{corol4} and Theorem \ref{theo2}, we have 
\begin{equation*}
\mathbf{P}\left( N_{t}\left( 1\right) =i\right) =\frac{\mu ^{t-1}\left( \mu
-1\right) ^{2}}{\left( \mu ^{t}-1\right) \left( \mu ^{t+1}-1\right) }\left[ 
\frac{\mu \left( \mu ^{t}-1\right) }{\mu ^{t+1}-1}\right] ^{i},
\end{equation*}
and 
\begin{equation*}
\mathbf{P}\left( \tau _{N_{t}\left( 1\right) ,1}^{\left( t\right) }\left(
1\right) \geq k\mid N_{t}\left( 1\right) >0\right) =\frac{\mu ^{t+1}-\mu ^{k}%
}{\mu ^{k}\left( \mu ^{t+1}-1\right) }.
\end{equation*}

This result follows from the linear-fractional model formulas 
\begin{eqnarray*}
\phi _{k}^{\prime }\left( z\right) &=&\frac{A_{k}D_{k}-B_{k}C_{k}}{\left(
C_{t}-D_{t}z\right) ^{2}} \\
\phi _{k}^{\prime }\left( \mathbf{P}\left( N_{t-k}\left( 1\right) =0\right)
\right) &=&\frac{A_{k}D_{k}-B_{k}C_{k}}{\left( C_{k}-D_{k}\frac{A_{t-k}}{%
C_{t-k}}\right) ^{2}}.\text{ }
\end{eqnarray*}

Depending on the value of $\mu $, these formulas yield the following results.

\begin{itemize}
\item  If $\mu <1$ (subcritical case), for $s\in \left\{ 2,...,t\right\} $%
\begin{equation*}
\mathbf{P}\left( t-\tau _{N_{t}\left( 1\right) ,1}^{\left( t\right) }\left(
1\right) =s\mid N_{t}\left( 1\right) >0\right) =\frac{\left( 1-\mu \right)
\mu ^{s}}{1-\mu ^{t+1}}\underset{t\rightarrow \infty }{\rightarrow }\left(
1-\mu \right) \mu ^{s},\text{ }s\geq 2
\end{equation*}
and if $s=1$, 
\begin{equation*}
\mathbf{P}\left( t-\tau _{N_{t}\left( 1\right) ,1}^{\left( t\right) }\left(
1\right) =1\mid N_{t}\left( 1\right) >0\right) \underset{t\rightarrow \infty 
}{\rightarrow }1-\mu ^{2}.
\end{equation*}
We conclude that $\tau _{N_{t}\left( 1\right) ,1}^{\left( t\right) }\left(
1\right) $ is of order $t$: coalescence time for the whole population
concentrates in the recent past generation $t$.\newline

\item  If $\mu >1$ (supercritical case), for $k\in \left\{ 0,...,t-2\right\} 
$%
\begin{equation*}
\mathbf{P}\left( \tau _{N_{t}\left( 1\right) ,1}^{\left( t\right) }\left(
1\right) =k\mid N_{t}\left( 1\right) >0\right) =\frac{\left( \mu -1\right)
\mu ^{t-k}}{\mu ^{t+1}-1}\underset{t\rightarrow \infty }{\rightarrow }\left(
\mu -1\right) \mu ^{-k-1},\text{ }k\geq 0,
\end{equation*}
and if $k=t-1$, 
\begin{equation*}
\mathbf{P}\left( \tau _{N_{t}\left( 1\right) ,1}^{\left( t\right) }\left(
1\right) =t-1\mid N_{t}\left( 1\right) >0\right) \underset{t\rightarrow
\infty }{\rightarrow }0.
\end{equation*}
$\tau _{N_{t}\left( 1\right) ,1}^{\left( t\right) }\left( 1\right) $ is of
order $1$, with geometric$\left( 1/\mu \right) $ law concentrated near the
remote past.\newline

\item  The previous results can be extended to the case $\mu =1$ (critical
case):

\begin{equation*}
\mathbf{P}\left( \tau _{i,1}^{\left( t\right) }\left( 1\right) \geq k\mid
N_{t}\left( 1\right) =i\right) =\left[ \frac{\left( t-k\right) \left(
t+1\right) }{t\left( t-k+1\right) }\right] ^{i-1}.
\end{equation*}
Letting $k=\left[ tx\right] $, for some $x\in \left( 0,1\right) $, we get 
\begin{equation*}
\mathbf{P}\left( \frac{\tau _{i,1}^{\left( t\right) }\left( 1\right) }{t}%
\geq x\mid N_{t}\left( 1\right) =i\right) \underset{t\rightarrow \infty }{%
\sim }1-\left( i-1\right) \frac{x}{1-x}\frac{1}{t}.
\end{equation*}

Moreover, 
\begin{equation*}
\mathbf{P}\left( \tau _{N_{t}\left( 1\right) ,1}^{\left( t\right) }\left(
1\right) \geq k\mid N_{t}\left( 1\right) >0\right) =\frac{t\left(
t-k+1\right) }{\left( t-1\right) \left( t+1\right) },
\end{equation*}
showing that 
\begin{equation*}
\mathbf{P}\left( \frac{\tau _{N_{t}\left( 1\right) ,1}^{\left( t\right)
}\left( 1\right) }{t}\geq x\mid N_{t}\left( 1\right) >0\right) \underset{%
t\rightarrow \infty }{\rightarrow }1-x,
\end{equation*}
the uniform distribution. This result appears in Theorem $2.1$ of \cite{A}.

In the critical case $\tau _{N_{t}\left( 1\right) ,1}^{\left( t\right)
}\left( 1\right) $ is of order $t$, with law concentrated in-between the
remote and recent past.
\end{itemize}

\subsubsection{\textbf{Triple versus binary one-step merging}}

The formulas for the binary and triple one-step merging read

\begin{eqnarray*}
\widehat{P}_{2,1}^{\left( t\right) }\left( 1\right) &=&2p^{2}\alpha
_{t}\left( \delta _{t}-\beta _{t}\right) \int_{0}^{1}dz\frac{1-z}{\left(
1-\delta _{t}z\right) ^{2}\left( 1-pz\right) ^{2}}, \\
\widehat{P}_{3,1}^{\left( t\right) }\left( 1\right) &=&3p^{3}\alpha
_{t}\left( \delta _{t}-\beta _{t}\right) \int_{0}^{1}dz\frac{\left(
1-z\right) ^{2}}{\left( 1-\delta _{t}z\right) ^{2}\left( 1-pz\right) ^{3}}.
\end{eqnarray*}
Simple calculations then lead to 
\begin{equation}
\begin{array}{l}
\frac{\widehat{P}_{3,1}^{\left( t\right) }\left( 1\right) }{\widehat{P}%
_{2,1}^{\left( t\right) }\left( 1\right) }\rightarrow \text{Constant}<1\text{
as }t\rightarrow \infty \text{ if }\mu <1, \\ 
\frac{\widehat{P}_{3,1}^{\left( t\right) }\left( 1\right) }{\widehat{P}%
_{2,1}^{\left( t\right) }\left( 1\right) }\sim \frac{3}{4}\frac{p\left(
2-p\right) }{\log t}\rightarrow 0\text{ as }t\rightarrow \infty \text{ if }%
\mu =1, \\ 
\frac{\widehat{P}_{3,1}^{\left( t\right) }\left( 1\right) }{\widehat{P}%
_{2,1}^{\left( t\right) }\left( 1\right) }\sim \frac{3}{4}\frac{p\left(
2-p\right) }{t\log \mu }\rightarrow 0\text{ as }t\rightarrow \infty \text{
if }\mu >1.\text{ }
\end{array}
\label{e16}
\end{equation}

In the subcritical case, the TMRCAs of 2 or 3 sampled individuals both are
of order $t$ (with a non-negligible probability to be $t-1$): the ratio $%
\widehat{P}_{3,1}^{\left( t\right) }\left( 1\right) /\widehat{P}%
_{2,1}^{\left( t\right) }\left( 1\right) $ goes to a constant limit which
can be explicitly computed. In the critical or in the supercritical case,
the TMRCA for 2 individuals has a low probability of being $t-1$, still
smaller when 3 individuals are being sampled: the ratio $\widehat{P}%
_{3,1}^{\left( t\right) }\left( 1\right) /\widehat{P}_{2,1}^{\left( t\right)
}\left( 1\right) $ goes to $0$ for large $t$; but at different speeds. In
the latter case, as $t$ gets large, only binary mergers will be seen.

These results follow from $\phi \left( z\right) =q/\left( 1-pz\right) $
leading to $\phi _{t}\left( z\right) =\alpha _{t}\frac{1-\beta _{t}z}{%
1-\delta _{t}z},$ $\phi _{t}^{\prime }\left( z\right) =\alpha _{t}\frac{%
\delta _{t}-\beta _{t}}{\left( 1-\delta _{t}z\right) ^{2}}$ and $\phi
^{\prime }\left( z\right) =qp/\left( 1-pz\right) ^{2}$, $\phi ^{\prime
\prime }\left( z\right) =2qp^{2}/\left( 1-pz\right) ^{3}$ and $\phi
^{^{\prime \prime \prime }}\left( z\right) =6qp^{3}/\left( 1-pz\right) ^{4}$.

These identities together with the Proposition \ref{prop1}, reading 
\begin{eqnarray*}
\widehat{P}_{i,1}^{\left( t\right) }\left( 1\right) &=&\frac{1}{\Gamma
\left( i\right) }\int_{0}^{1}dz\left( 1-z\right) ^{i-1}\left[ \phi
_{t-1}\right] ^{\left( 1\right) }\left( \phi \left( z\right) \right) \phi
^{\left( i\right) }\left( z\right) \\
&=&\frac{1}{\Gamma \left( i\right) }\int_{0}^{1}dz\left( 1-z\right) ^{i-1}%
\frac{\phi _{t}^{\left( 1\right) }\left( z\right) }{\phi \left( z\right) }%
\phi ^{\left( i\right) }\left( z\right) ,
\end{eqnarray*}
give the result.

\subsection{Sibuya branching mechanism with infinite mean $\mu $}

Let $\alpha \in \left( 0,1\right) $ and $\lambda \in \left( 0,1\right] $.
With $\left[ a\right] _{i}=a\left( a+1\right) ...\left( a+i-1\right) $ the $%
i-$th rising factorial of $a$, consider a branching mechanism with
probability mass 
\begin{equation*}
\pi _{0}:=\mathbf{P}\left( M=0\right) =1-\lambda \text{ and}
\end{equation*}
\begin{equation}
\pi _{m}:=\mathbf{P}\left( M=m\right) =\lambda \left( -1\right) ^{m-1}\binom{%
\alpha }{m}=\alpha \lambda \frac{\left[ 1-\alpha \right] _{m-1}}{m!}\text{, }%
m\geq 1.  \label{e17}
\end{equation}
The pgf of $M$ is $\phi \left( z\right) =1-\lambda \left( 1-z\right)
^{\alpha }$, leading to\footnote{%
This in accordance with the arguments developed in the proof of Proposition $%
2$ page $3764$ of \cite{A2}.}: $\phi _{t}\left( z\right) =1-\lambda
_{t}\left( 1-z\right) $, with $\lambda _{t}=\lambda ^{\left( 1-\alpha
^{t}\right) /\left( 1-\alpha \right) }$ and $\alpha _{t}=\alpha ^{t}$. We
thus have $\mathbf{P}\left( N_{t}\left( 1\right) =0\right) =1-\lambda _{t}$
and

\begin{equation*}
\mathbf{P}\left( N_{t}\left( 1\right) =n\right) =\lambda _{t}\left(
-1\right) ^{n-1}\binom{\alpha _{t}}{n}=\alpha _{t}\lambda _{t}\frac{\left[
1-\alpha _{t}\right] _{n-1}}{n!}\text{, }n\geq 1.
\end{equation*}
When $\lambda =1$, $M$ can be seen as the first epoch of a success in a
Bernoulli trial when the probability of success is inversely proportional to
the number of the trial, see e.g. \cite{Hui}.\newline

\emph{Remark 4:} One can easily check part of the Corollary $4$, while
observing: 
\begin{eqnarray*}
\mathbf{P}\left( N_{t}\left( 1\right) \geq 2\right) &=&\lambda _{t}\left(
1-\alpha _{t}\right) , \\
\widehat{P}_{2,1}^{\left( t\right) }\left( n_{0}\right) &=&\lambda
_{t}\left( 1-\alpha \right) , \\
\widehat{P}_{2,2}^{\left( t\right) }\left( n_{0}\right) &=&\lambda
_{t}\left( \alpha -\alpha _{t}\right) .\text{ }\Diamond
\end{eqnarray*}

\subsubsection{\textbf{Coalescence for pairs of tips: tail probability}}

The TMRCA distribution, assuming it is finite, is $\mathbf{P}\left( t-\tau
_{2,1}^{\left( t\right) }\left( n_{0}\right) =l\mid \tau _{2,1}^{\left(
t\right) }\left( n_{0}\right) <\infty \right) =\frac{\alpha ^{l-1}\left(
1-\alpha \right) }{1-\alpha ^{t}}$, $l=1,...,t.$ Thus 
\begin{equation}
t-\tau _{2,1}^{\left( t\right) }\left( n_{0}\right) \mid \tau _{2,1}^{\left(
t\right) }\left( n_{0}\right) <\infty \overset{d}{\underset{t\rightarrow
\infty }{\rightarrow }}\text{geom}\left( \alpha \right) ,  \label{e18}
\end{equation}
a geometric distribution with success probability $\alpha $. Note 
\begin{equation}
\begin{array}{l}
\mathbf{P}\left( t-\tau _{2,1}^{\left( t\right) }\left( n_{0}\right)
=l\right) =\left( 1-\left( 1-\lambda _{t}\right) ^{n_{0}}\right) \alpha
^{l-1}\left( 1-\alpha \right) ,\text{ }l=1,...,t \\ 
\rightarrow \left( 1-\left( 1-\lambda ^{1/\left( 1-\alpha \right) }\right)
^{n_{0}}\right) \alpha ^{l-1}\left( 1-\alpha \right) ,\text{ }l\geq 1\text{
as }t\rightarrow \infty .
\end{array}
\label{e19}
\end{equation}
If $\lambda \neq 1$, the limit law is defective in that its total mass is $%
1-\left( 1-\lambda ^{1/\left( 1-\alpha \right) }\right) ^{n_{0}}<1$, which
is the probability of non-extinction of the branching process with branching
mechanism $\phi \left( z\right) =1-\lambda \left( 1-z\right) ^{\alpha }$ and 
$n_{0}$ founders.

The coalescence time for a randomly chosen pair is geometrically
concentrated near the recent past generation $t$.

This result follows from the Corollary \ref{corol2}: 
\begin{eqnarray*}
&&\mathbf{P}\left( \infty >\tau _{2,1}^{\left( t\right) }\left( n_{0}\right)
\geq k\right) \\
&=&n_{0}\int_{0}^{1}dz\left( 1-z\right) \frac{\phi _{t-k}^{\prime \prime
}\left( z\right) }{\phi _{t-k}^{\prime }\left( z\right) }\phi _{t}^{\prime
}\left( z\right) \phi _{t}\left( z\right) ^{n_{0}-1} \\
&=&n_{0}\alpha _{t}\lambda _{t}\left( 1-\alpha _{t-k}\right)
\int_{0}^{1}du\cdot u^{\alpha _{t}-1}\left( 1-\lambda _{t}u^{\alpha
_{t}}\right) ^{n_{0}-1} \\
&=&n_{0}\left( 1-\alpha _{t-k}\right) \int_{0}^{\lambda _{t}}dv\cdot \left(
1-v\right) ^{n_{0}-1} \\
&=&\left( 1-\alpha _{t-k}\right) \left( 1-\left( 1-\lambda _{t}\right)
^{n_{0}}\right) ,
\end{eqnarray*}
and 
\begin{eqnarray*}
\mathbf{P}\left( \tau _{2,1}^{\left( t\right) }\left( n_{0}\right) =k\right)
&=&\left( 1-\left( 1-\lambda _{t}\right) ^{n_{0}}\right) \alpha
^{t-k-1}\left( 1-\alpha \right) \text{, }k=0,...,t-1 \\
\mathbf{P}\left( \tau _{2,1}^{\left( t\right) }\left( n_{0}\right) <\infty
\right) &=&\sum_{k=0}^{t-1}\mathbf{P}\left( \tau _{2,1}^{\left( t\right)
}\left( n_{0}\right) =k\right) =\left( 1-\left( 1-\lambda _{t}\right)
^{n_{0}}\right) \left( 1-\alpha ^{t}\right) \\
\mathbf{P}\left( \tau _{2,1}^{\left( t\right) }\left( n_{0}\right) =k\mid
\tau _{2,1}^{\left( t\right) }\left( n_{0}\right) <\infty \right) &=&\frac{%
\alpha ^{t-k-1}\left( 1-\alpha \right) }{1-\alpha ^{t}}.
\end{eqnarray*}

\subsubsection{\textbf{Existence of a common ancestor}}

$\tau _{2,1}^{\left( t\right) }\left( n_{0}\right) $ can take the value $%
+\infty $ for two reasons: one is because the two sampled particles do not
belong to the descendance of the same founder, the other one being because
there could be less than two particles to be sampled at generation $t $, the
latter event occurring with probability 
\begin{eqnarray*}
\mathbf{P}\left( N_{t}\left( n_{0}\right) <2\right) &=&\phi _{t}\left(
0\right) ^{n_{0}}+n_{0}\phi _{t}^{\prime }\left( 0\right) \phi _{t}\left(
0\right) ^{n_{0}-1} \\
&=&\left( 1-\lambda _{t}\right) ^{n_{0}}+n_{0}\alpha _{t}\lambda _{t}\left(
1-\lambda _{t}\right) ^{n_{0}-1} \\
&=&\left( 1-\lambda _{t}\right) ^{n_{0}-1}\left( 1+\lambda _{t}\left(
n_{0}\alpha ^{t}-1\right) \right) .
\end{eqnarray*}
To be complete, one needs to compute the probability that $\tau
_{2,1}^{\left( t\right) }\left( n_{0}\right) $ takes the value $+\infty $
resulting only from the two sampled particles not belonging to the
descendance of the same founder. We find 
\begin{eqnarray*}
&&\mathbf{P}\left( N_{t}\left( n_{0}\right) \geq 2,\tau _{2,1}^{\left(
t\right) }\left( n_{0}\right) =\infty \right) \\
&=&n_{0}\left( n_{0}-1\right) \int_{0}^{1}dz\left( 1-z\right) \phi
_{t}^{\prime }\left( z\right) ^{2}\phi _{t}\left( z\right) ^{n_{0}-2} \\
&=&\mathbf{P}\left( N_{t}\left( n_{0}\right) \geq 2\right) -\mathbf{P}\left(
\tau _{2,1}^{\left( t\right) }\left( n_{0}\right) <\infty \right) \\
&=&1-\left( 1-\lambda _{t}\right) ^{n_{0}}-n_{0}\alpha _{t}\lambda
_{t}\left( 1-\lambda _{t}\right) ^{n_{0}-1}-\left( 1-\left( 1-\lambda
_{t}\right) ^{n_{0}}\right) \left( 1-\alpha ^{t}\right) \\
&=&\alpha ^{t}\left( 1-\left( 1-\lambda _{t}\right) ^{n_{0}-1}\left(
1+\left( n_{0}-1\right) \lambda _{t}\right) \right) ,
\end{eqnarray*}
to be compared with the full probability 
\begin{equation*}
\mathbf{P}\left( \tau _{2,1}^{\left( t\right) }\left( n_{0}\right) =\infty
\right) =\alpha ^{t}+\left( 1-\lambda _{t}\right) ^{n_{0}}\left( 1-\alpha
^{t}\right) .
\end{equation*}
Note that if $\lambda =1$ (whence $\lambda _{t}=1$) and $n_{0}\geq 2$, $%
N_{t}\left( n_{0}\right) \geq n_{0}\geq 2$ and then $\mathbf{P}\left( \tau
_{2,1}^{\left( t\right) }\left( n_{0}\right) =\infty \right) =\alpha ^{t}$,
independently of the true value of $n_{0}$: the two sampled particles belong
to the descendance of the same founder with probability close to $1$ (the
founder with largest family size being dominant) as $t$ increases. In this
case, 
\begin{equation}
\begin{array}{l}
\mathbf{P}\left( t-\tau _{2,1}^{\left( t\right) }\left( n_{0}\right)
=l\right) =\alpha ^{l-1}\left( 1-\alpha \right) \text{, }l=1,...,t\text{,
with} \\ 
t-\tau _{2,1}^{\left( t\right) }\left( n_{0}\right) \rightarrow \text{geom}%
\left( \alpha \right) \text{ as }t\rightarrow \infty \text{, in distribution.%
}
\end{array}
\label{e20}
\end{equation}
As suggested in \cite{A2}, for such rapidly growing populations with $%
\lambda =1$, coalescence occurs in the very recent past close to $t.$

Due to heavy-tailedness of the offspring distribution, it is very likely
that two randomly chosen individuals belong to the same family (the
descendance at $t$ of one of the $n_{0}$ founders) and within this family,
it is very likely that the two randomly chosen individuals belong to the
direct offspring of the largest family size at $t$ of this founder. The
probability that this does not happen is small, of order $\alpha ^{t}$, if $%
t $ is large.

\subsubsection{\textbf{TMRCA for the whole population currently alive}}

Let us finally consider the problem of the TMRCA for the whole population
alive at $t$. We have 
\begin{equation*}
\mathbf{P}\left( \tau _{i,1}^{\left( t\right) }\left( 1\right) \geq k\mid
N_{t}\left( 1\right) =i\right) =\frac{\left[ 1-\alpha ^{t-k}\right] _{i-1}}{%
\left[ 1-\alpha ^{t}\right] _{i-1}}.
\end{equation*}
If $\alpha \ll 1_{,}$ $\left[ 1-\alpha \right] _{i-1}\sim i!\left( 1-\alpha
H_{i}\right) $ where $H_{i}$ is the Harmonic number. Thus, if $\alpha
^{t-k}\ll 1$, with $s=t-k$, as $t\rightarrow \infty $, 
\begin{equation}
\begin{array}{l}
\mathbf{P}\left( \tau _{i,1}^{\left( t\right) }\left( 1\right) \geq k\mid
N_{t}\left( 1\right) =i\right) \sim 1-\alpha ^{t}\left( \alpha
^{-k}-1\right) H_{i} \\ 
\mathbf{P}\left( t-\tau _{i,1}^{\left( t\right) }\left( 1\right) =s\mid
N_{t}\left( 1\right) =i\right) \sim \alpha ^{s-1}\left( 1-\alpha \right)
H_{i},
\end{array}
\label{e21}
\end{equation}
showing that $\tau _{i,1}^{\left( t\right) }\left( 1\right) $ is `close to $%
t $'.

Furthermore, for all $t\geq 1$, 
\begin{equation}
\begin{array}{l}
\mathbf{P}\left( \tau _{N_{t}\left( 1\right) ,1}^{\left( t\right) }\left(
1\right) =k\mid N_{t}\left( 1\right) >0\right) =\left( 1-\alpha \right)
\alpha ^{k}\text{, }k=0,...,t-2 \\ 
\mathbf{P}\left( \tau _{N_{t}\left( 1\right) ,1}^{\left( t\right) }\left(
1\right) =t-1\mid N_{t}\left( 1\right) >0\right) =\alpha ^{t-1},
\end{array}
\label{e22}
\end{equation}
approaching the geometric$\left( \alpha \right) $ distribution as $t$ gets
large. The order of magnitude of $\tau _{N_{t}\left( 1\right) ,1}^{\left(
t\right) }\left( 1\right) $ is $1$ and $\tau _{N_{t}\left( 1\right)
,1}^{\left( t\right) }\left( 1\right) $ is geometrically concentrated near
the root of the tree (in the remote past).

This result follows from the following facts: the probability conditioned on 
$N_{t}(1)=i$ is given by Corollary \ref{corol4}, which reads 
\begin{equation*}
\mathbf{P}\left( \tau _{i,1}^{\left( t\right) }\left( 1\right) \geq k\mid
N_{t}\left( 1\right) =i\right) =\frac{\left( \phi _{k}\right) ^{\left(
1\right) }\left( \mathbf{P}\left( N_{t-k}\left( 1\right) =0\right) \right) 
\mathbf{P}\left( N_{t-k}\left( 1\right) =i\right) }{\mathbf{P}\left(
N_{t}\left( 1\right) =i\right) }.
\end{equation*}
The model yields the formulas $\phi _{t}\left( z\right) =1-\lambda
_{t}\left( 1-z\right) ^{\alpha _{t}},$ with $\lambda _{t}=\lambda ^{\left(
1-\alpha ^{t}\right) /\left( 1-\alpha \right) }$ and $\alpha _{t}=\alpha
^{t} $, 
\begin{equation*}
\mathbf{P}\left( N_{t}\left( 1\right) =0\right) =1-\lambda _{t}\text{ and}
\end{equation*}
\begin{equation}
\mathbf{P}\left( N_{t}\left( 1\right) =i\right) =\lambda _{t}\left(
-1\right) ^{i-1}\binom{\alpha _{t}}{i}=\alpha _{t}\lambda _{t}\frac{\left[
1-\alpha _{t}\right] _{i-1}}{i!}\text{, }k\geq 1.  \label{sprob}
\end{equation}

\begin{eqnarray*}
\phi _{k}^{\prime }\left( z\right) &=&\alpha _{k}\lambda _{k}\left(
1-z\right) ^{\alpha _{k}-1}, \\
\phi _{k}^{\prime }\left( \mathbf{P}\left( N_{t-k}\left( 1\right) =0\right)
\right) &=&\alpha _{k}\lambda _{k}\lambda _{t-k}^{\alpha _{k}-1}.
\end{eqnarray*}
The probability conditioned on $N_{t}(1)=i$ follows from Theorem \ref{theo2}%
, 
\begin{eqnarray*}
\mathbf{P}\left( \tau _{N_{t}\left( 1\right) ,1}^{\left( t\right) }\left(
1\right) \geq k\mid N_{t}\left( 1\right) >0\right) &=&\frac{\mathbf{P}\left(
N_{t-k}\left( 1\right) >0\right) }{\mathbf{P}\left( N_{t}\left( 1\right)
>0\right) }\left( \phi _{k}\right) ^{\left( 1\right) }\left( \mathbf{P}%
\left( N_{t-k}\left( 1\right) =0\right) \right) \\
&=&\frac{\lambda _{t-k}}{\lambda _{t}}\alpha _{k}\lambda _{k}\lambda
_{t-k}^{\alpha _{k}-1}=\alpha ^{k},\text{ }k=0,...,t-1.\text{ }
\end{eqnarray*}

\subsubsection{\textbf{Triple versus binary one-step mergers}}

The formula for a one-step merging reads $\widehat{P}_{i,1}^{\left( t\right)
}\left( 1\right) =\frac{\Gamma \left( i-\alpha \right) }{\Gamma \left(
1-\alpha \right) \Gamma \left( i\right) }\lambda ^{\left( 1-\alpha
^{t}\right) /\left( 1-\alpha \right) }.$ Thus, $\widehat{P}_{2,1}^{\left(
t\right) }\left( 1\right) =\left( 1-\alpha \right) \lambda ^{\left( 1-\alpha
^{t}\right) /\left( 1-\alpha \right) }$ and $\widehat{P}_{3,1}^{\left(
t\right) }\left( 1\right) =\left( 1-\alpha \right) \left( 1-\alpha /2\right)
\lambda ^{\left( 1-\alpha ^{t}\right) /\left( 1-\alpha \right) }$, leading
to 
\begin{equation}
\frac{\widehat{P}_{3,1}^{\left( t\right) }\left( 1\right) }{\widehat{P}%
_{2,1}^{\left( t\right) }\left( 1\right) }\rightarrow \left( 1-\alpha
/2\right) \text{ as }t\rightarrow \infty .  \label{e23}
\end{equation}
For very rapidly growing populations, the ratio $\widehat{P}_{3,1}^{\left(
t\right) }\left( 1\right) /\widehat{P}_{2,1}^{\left( t\right) }\left(
1\right) $ goes to a limit for large $t$.\newline

This result follows from the following derivations:

$\phi \left( z\right) =1-\lambda \left( 1-z\right) ^{\alpha }$, leading to $%
\phi _{t}\left( z\right) =1-\lambda _{t}\left( 1-z\right) ^{\alpha _{t}}$
with $\lambda _{t}=\lambda ^{\left( 1-\alpha ^{t}\right) /\left( 1-\alpha
\right) }$ and $\alpha _{t}=\alpha ^{t},$ we have

$\phi _{t-1}^{\prime }\left( z\right) =\alpha _{t-1}\lambda _{t-1}\left(
1-z\right) ^{\alpha _{t-1}-1}$ and $\phi _{t-1}^{\prime }\left( \phi \left(
z\right) \right) =\alpha _{t-1}\lambda _{t-1}\left( \lambda \left(
1-z\right) ^{\alpha }\right) ^{\alpha _{t-1}-1},$

$\phi ^{\prime }\left( z\right) =\alpha \lambda \left( 1-z\right) ^{\alpha
-1}$, $\phi ^{\prime \prime }\left( z\right) =\alpha \left( 1-\alpha \right)
\lambda \left( 1-z\right) ^{\alpha -2}$ and

$\phi ^{^{\prime \prime \prime }}\left( z\right) =\alpha \left( 1-\alpha
\right) \left( 2-\alpha \right) \lambda \left( 1-z\right) ^{\alpha -3}$, and
finally

\begin{equation*}
\widehat{P}_{i,1}^{\left( t\right) }\left( 1\right) =\frac{1}{\Gamma \left(
i\right) }\int_{0}^{1}dz\left( 1-z\right) ^{i-1}\left[ \phi _{t-1}\right]
^{\left( 1\right) }\left( \phi \left( z\right) \right) \phi ^{\left(
i\right) }\left( z\right) ,\text{ }i\geq 2.
\end{equation*}
\newline

\textbf{Acknowledgments:}

T. Huillet acknowledges partial support from the ``Chaire \textit{%
Mod\'{e}lisation math\'{e}matique et biodiversit\'{e}''.} N. Grosjean and T.
Huillet also acknowledge support from the labex MME-DII Center of Excellence
(\textit{Mod\`{e}les math\'{e}matiques et \'{e}conomiques de la dynamique,
de l'incertitude et des interactions}, ANR-11-LABX-0023-01 project).

\end{document}